\documentclass[reqno]{amsart}[standalone]
\usepackage[notocbasic]{nomencl}
\usepackage{diagbox}
\makenomenclature
\usepackage{algorithm2e}  %
\RestyleAlgo{boxruled}
\makeatletter
\newcommand{\mycaption}{%
\ifx \@captype \@undefined \@latex@error {\noexpand \caption outside float}\@ehd \expandafter \@gobble \else \refstepcounter \@captype \expandafter \@firstofone \fi {\@dblarg {\@caption \@captype }}%
}%
\makeatother
\usepackage{pifont}
\newcommand{\cmark}{\ding{51}}
\newcommand{\xmark}{\ding{55}}
\usepackage[table]{xcolor}
\usepackage{tikz}
\usepackage{multirow}
\usetikzlibrary{positioning, shapes.geometric}
\usepackage{units}
\usepackage{amssymb}
\usepackage{amsmath}
\usepackage{mathrsfs}
\usepackage[foot]{amsaddr}
\usepackage{amsfonts}
\usepackage{bm}
\usepackage{amsmath}
\usepackage{graphicx}
\usepackage{booktabs}
\usepackage{caption}
\usepackage{subcaption}
\usepackage{mathbbol}
\usepackage{hyperref}
\usepackage{lipsum}
\usepackage{mathtools}

\newcommand{\newinf}{\mathop{\mathrm{inf}\vphantom{\mathrm{sup}}}}
\usepackage{hhline}
\usepackage{pgf}

\DeclareMathOperator*{\argmin}{arg\,min}
\graphicspath{{Images/}}
\usepackage[top=1in, bottom=1.25in, left=1.25in, right=1.25in]{geometry}
\usepackage [autostyle, english = american]{csquotes}
\MakeOuterQuote{"}
\newcommand{\ie}{{\it i.e.}}
\definecolor{vargreen}{rgb}{0.0, 0.5, 0.0}
\definecolor{vared}{rgb}{0.9, 0.0, 0.0}
\definecolor{low}{rgb}{0.37, 0.79, 0.38}
\definecolor{high}{rgb}{0.27, 0., 0.33}

\usepackage{tikz}
\usetikzlibrary{babel}
\usepackage{etoolbox} %
\usepackage{listofitems} %
\usepackage{pgfmath}
\usepackage{todonotes} %
\newtheorem{remark}{Remark}
\begin{document}
\title[OT-inspired framework for DL in ROMs]{Optimal Transport-inspired Deep Learning Framework for Slow-Decaying Kolmogorov n-width Problems: Exploiting Sinkhorn Loss and Wasserstein Kernel}
\author{Moaad Khamlich$^{1}$}
\author{Federico Pichi$^{2}$}
\author{Gianluigi Rozza$^{1}$}
\address{$^1$ mathLab, Mathematics Area, SISSA, via Bonomea 265, I-34136 Trieste, Italy}
\address{$^2$ Chair of Computational Mathematics and Simulation Science, \'Ecole Polytechnique F\'ed\'erale de Lausanne, 1015 Lausanne, Switzerland}
\maketitle
\begin{abstract}
Reduced order models (ROMs) are widely used in scientific computing to tackle high-dimensional systems. However, traditional ROM methods may only partially capture the intrinsic geometric characteristics of the data. These characteristics encompass the underlying structure, relationships, and essential features crucial for accurate modeling.

To overcome this limitation, we propose a novel ROM framework that integrates optimal transport (OT) theory and neural network-based methods.
Specifically, we investigate the Kernel Proper Orthogonal Decomposition (kPOD) method exploiting the Wasserstein distance as the custom kernel, and we efficiently train the resulting neural network (NN) employing the Sinkhorn algorithm. By leveraging an OT-based nonlinear reduction, the presented framework can capture the geometric structure of the data, which is crucial for accurate learning of the reduced solution manifold.
When compared with traditional metrics such as mean squared error or cross-entropy, exploiting the Sinkhorn divergence as the loss function enhances stability during training, robustness against overfitting and noise, and accelerates convergence.

To showcase the approach's effectiveness, we conduct experiments on a set of challenging test cases exhibiting a slow decay of the Kolmogorov n-width. The results show that our framework outperforms traditional ROM methods in terms of accuracy and computational efficiency.

\end{abstract}
\section{Introduction}
\label{sec:intro}
Reduced order modeling has emerged as a powerful technique to speed up the
investigation of complex systems described by partial differential equations
(PDEs). In fact, the discretization of generic nonlinear PDEs often yields
high-dimensional solution vectors, meaning that the description of their
behavior is well-captured when employing a consistent number of degrees of
freedom. This can have a substantial impact on the computational cost of
analyzing the physical phenomena.

Fortunately, many practical problems exhibit a low dimensional set of dominant
patterns that characterize the solution's variation with respect to the
parameters, indicating the presence of redundant variables in the description of
the phenomena. This motivates the use of dimensionality reduction techniques to
compress the original dataset while retaining its essential properties. By
employing these techniques, it becomes possible to create low-dimensional
surrogate models, known as reduced order models (ROMs)\nomenclature{ROM}{Reduced Order Model} \cite{benner2017model,benner2017, Grepl2007,schilders,vol1,vol2}, that capture the essential features of the original
high-dimensional full order model (FOM)\nomenclature{FOM}{Full Order Model},
while reducing the computational cost of the simulation.

Specifically, dimensionality reduction techniques can be classified into linear
and nonlinear methods, depending on the nature of the map from the reduced order to
the full order space.

Linear approaches, such as the proper orthogonal decomposition (POD) and the Greedy algorithm, measure the
approximability of the solution within a linear trial subspace in terms of the Kolmogorov $n$-width.
The Kolmogorov $n$-width, denoted by $d_n(\mathcal{S}, |\cdot|)$, is a measure of the best possible approximation of a set $\mathcal{S}$ by an $n$-dimensional subspace $\mathcal{V}$, where $|\cdot|$ represents a norm. It quantifies how well the elements in $\mathcal{S}$ can be approximated using elements from $\mathcal{V}$, and it is defined as:

\begin{equation}
d_n(\mathcal{S}, |\cdot|) \stackrel{\text { def }}{=} \newinf_{\substack{\mathcal{V}\ \text{s.t.} \\ \dim(\mathcal{V}) = n}} \sup_{u \in \mathcal{S}} \newinf_{v \in \mathcal{V}} |u - v|.
\end{equation}
Handling problems with slow decay of the Kolmogorov $n$-width poses a challenge in the design of efficient ROMs, as increasing the
dimensionality of the reduced model only allows to achieve a slight gain in
accuracy.

In particular, despite its clarity and interpretability, POD suffers from
certain drawbacks, such as \textit{(i)} information loss affecting the modes,  and \textit{(ii)} the challenge of representing
moving discontinuities or advection-dominated effects. These issues are
particularly prominent in nonlinear problems, where complex interactions
usually introduce features that a linear POD basis struggles to reconstruct.
Moreover, in the nonlinear case, hyper-reduction procedures are often necessary
to make the reduced equations independent of the number of degrees of freedom
of the FOM. In such cases, nonlinear reduction techniques are more capable of
spanning the manifold associated with a parameterized PDE, which is
why they have become increasingly popular in the ROM community, dealing with
the effort of moving towards complex real-life simulations.

One notable advancement in nonlinear reduction techniques involves the
utilization of deep learning (DL) frameworks \nomenclature{DL}{Deep Learning}, which efficiently extract relevant features from high-dimensional
data \cite{goodfellow,machine-learning-book}. They overcome the constraints of traditional projection-based approaches by compressing the original dataset into a small latent space, albeit with
an increased training cost during the \textit{offline} stage.

In this work, we propose a novel DL-based ROM framework that leverages the
powerful tools of Optimal Transport (OT) theory and Kernel Proper Orthogonal
Decomposition (kPOD). OT provides a principled way to measure the distance
between probability distributions, which is a fundamental problem in many areas of science and engineering \cite{santambrogio2015optimal,peyre2019computational}. Our framework exploits the OT-based distance metric, known as the Wasserstein distance, to construct a custom kernel that captures the underlying hidden features of the data. The dimensionality reduction is carried out through a nonlinear mapping using kPOD \cite{kpca}, exploiting the kernel trick to implicitly employ a high-dimensional space known as the \textit{feature space}.

To train our models, \ie, to recover the forward and backward mappings between the full order solution and
latent space, we exploit an autoencoder
architecture, where the encoding layer is forced to learn the reduced
representation. The training of the proposed architecture is performed using
the Sinkhorn algorithm, a computationally efficient method for solving the OT
problem. Specifically, we use the Sinkhorn loss as a regularization term to
encourage the learned representations to be invariant to the input domain. Our
numerical investigation shows promising results, achieving state-of-the-art
performance when compared to existing DL-based ROM methods.

The remainder of this paper is organized as follows. Section
\ref{sec:framework_overview} introduces the proposed framework. Section
\ref{sec:related} provides a brief overview of related works focusing on DL
and OT-based ROMs. Sections \ref{sec:ot} and \ref{sec:kpod} give a general
introduction to OT theory and the kPOD procedure, respectively. In Section
\ref{sec:complete_framework}, we describe the main components of our framework
in light of the theory presented in the previous sections. Finally, Section
\ref{sec:numerical_results} presents numerical experiments on benchmark problems to evaluate the effectiveness of the proposed strategy.
\section{Framework Overview}
\label{sec:framework_overview}
In this section, we introduce the framework for reduced-order modeling in the
context of solving a generic parametrized PDE. We start by considering the
abstract PDE defined on a domain $\Omega \subset \mathbb{R}^n$:
\begin{equation}
    \label{eq:abs}
    \mathcal{L}(u(\boldsymbol{\mu}), \boldsymbol{\mu}) = 0,
\end{equation}
where $u(\boldsymbol{\mu}) \in \mathcal{V} $ is the unknown function belonging
to a suitable function space, and $\boldsymbol{\mu} \in \mathcal{P} \subset
\mathbb{R}^{p}$ represents the parameter vector. The operator $\mathcal{L}$
incorporates the differential operators, including boundary conditions and forcing terms, that govern the physical behavior of $u(\boldsymbol{\mu})$.
Discretization methods such as
finite differences, finite elements, or finite volumes are commonly used for
the numerical approximation of the PDE. The resulting discretized solution is represented by a vector
$\boldsymbol{u}_h(\boldsymbol{\mu}) \in \mathbb{R}^{N_h}$, where $N_h$ is the
number of degrees of freedom.
In particular, there is often a need to solve
\eqref{eq:abs} for different parameter values, \ie, in a many-query context. The investigation of the PDE can be pursued by exploring the so-called
solution manifold\footnote{While it is common to refer to $\mathcal{S}_h$
as a \textit{"manifold"}, it may not satisfy the definition of a manifold in
the context of differential geometry.}:

\begin{equation} \mathcal{S}_h = \left\{\boldsymbol{u}_h(\boldsymbol{\mu}) \mid
        \boldsymbol{\mu} \in \mathcal{P}\right\}\footnote{The model order
        reduction framework extends to time-dependent problems, in which
        $\boldsymbol{u}_{h} = \boldsymbol{u}_{h}(t, \boldsymbol{\mu})$.  Considering a data-driven approach, to simplify the notation, we incorporate the time variable as an additional parameter.}.  \end{equation}

As explained in Section \ref{sec:intro}, the discretization of Problem
\eqref{eq:abs} involves high-dimensional solution vectors, requiring many degrees of freedom ($N_h \gg 1$) to describe its behavior
accurately. This poses challenges in terms of computational cost for obtaining and analyzing
the data. Therefore, dimensionality reduction techniques are employed to
compress the information while preserving its main properties.

\begin{figure}[h] \centering \includegraphics[width=1\linewidth]{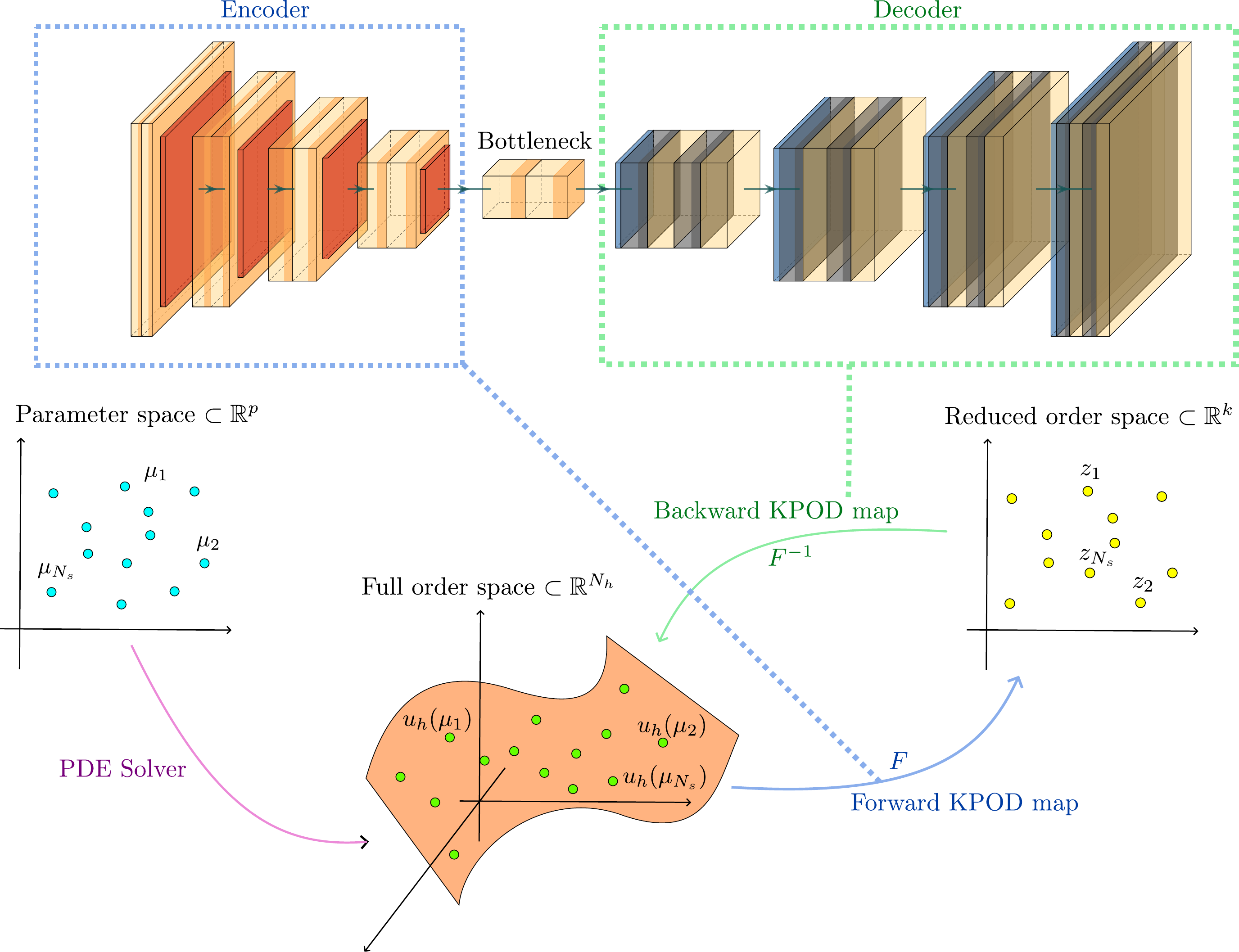}
    \caption{Schematic illustration of the different steps of the proposed
    methodology} \label{fig:scheme} \end{figure}

Since many reduction techniques rely on matrix-structured operations,
let us consider $N_s$ sample points in $\mathcal{P}$, and organize the
corresponding PDE solutions column-wise in the \textit{snapshot matrix}:
\begin{equation}
    \label{eq:snap}
    \mathbf{S} = \left[\boldsymbol{u}_h(\boldsymbol{\mu}_1) \mid
\boldsymbol{u}_h(\boldsymbol{\mu}_2) \mid \ldots \mid
\boldsymbol{u}_h(\boldsymbol{\mu}_{N_s})\right] \in \mathbb{R}^{N_h \times
N_s}.
\end{equation}
In particular, our dimensionality reduction framework, depicted in Figure \ref{fig:scheme}, transforms the dataset
$\mathbf{S}$ into a new dataset $\mathbf{Z} \in \mathbb{R}^{k\times N{s}}$,
while preserving the geometric structure of the data. The
columns of $\mathbf{Z}$ can be interpreted as images of the original FOM
snapshots $\{\boldsymbol{u}_{h}(\boldsymbol{\mu}_{j})\}_{j=1}^{N_{s}}$ through
an appropriate projection map $F: \boldsymbol{u}_{h} \in \mathbb{R}^{N_{h}}\rightarrow \boldsymbol{z} \in \mathbb{R}^{k}$, which we will refer to as the
\textit{forward map}. To reconstruct the full-order solution from the reduced
representation a \textit{backward map} $F^{-1}:\mathbb{R}^k \rightarrow
\mathbb{R}^{N_h}$ is required.

For example, the POD technique based on the singular value decomposition (SVD)
of the matrix $\mathbf{S}$, performs dimensionality reduction using the forward
map given by $F(\boldsymbol{u}_{h})= \mathbf{V}^{\top} \mathbf{u}_h$, where the orthonormal
matrix $\mathbf{V}$ spans the reduced basis. The advantage
of POD is that the backward map is also linear, and can be expressed as
$F^{-1}(\boldsymbol{z}) = \mathbf{V} \boldsymbol{z}$.

The standard ROM setting, within the context of solving the parametric Problem
\eqref{eq:abs}, aims at constructing a map from the parameter space to the reduced space $\mathscr{F}: \boldsymbol{\mu} \in \mathcal{P} \rightarrow \mathbf{u}_{r}(\boldsymbol{\mu}) \in \mathbb{R}^{N{r}}$. By utilizing the
compressed information, which contains $N_{s}$ known points of
$\mathscr{F}$, various projection/interpolation/regression techniques
can be employed to evaluate the reduced coefficients at other parameter values.

To overcome the limitations of linear methods for slow decaying
Kolmogorov $n$-width problems, our novel ROM framework is based on Optimal Transport (OT) theory (see Section \ref{sec:ot}) and a nonlinear extension of POD, known as Kernel POD (kPOD) \cite{kpca} (see Section \ref{sec:kpod}).

The latter consists of applying POD in a high-dimensional implicit space with higher expressive capabilities, by introducing a suitable bivariate symmetric form $\kappa:
\mathbb{R}^{N_{h}}\times\mathbb{R}^{N_{h}}\mapsto\mathbb{R}$, known as the
\textit{kernel function}. This kernel function is employed to define the kPOD
forward map $F$, which specifies a nonlinear projection of the FOM solutions
onto the kPOD principal components. It should be noted that, unlike the POD
procedure, here the projection is performed without having an explicit
expression of the principal components, \ie, the high-dimensional feature space
is only exploited implicitly through the kernel function.
When the kernel function is chosen as the standard inner product $\langle
\cdot, \cdot \rangle$, the kPOD framework coincides with classical POD.
However, to introduce nonlinearity into the projection process, and capture complex underlying features in the data, we employ a kernel function based on the \textit{Wasserstein distance}.  This distance measure is an important tool from OT theory that can be used to compare probability measures supported on high-dimensional spaces.

Given the nonlinear nature of the reduction performed by kPOD, constructing the
backward map $F^{-1}$ becomes a nontrivial task. To address this challenge, we
employed a deep learning approach. Specifically, we used an autoencoder
architecture, which aims at reconstructing the identity map $\mathcal{I}$, to find an approximation $\tilde{\boldsymbol{u}}_{h}$ of $\boldsymbol{{u}}_{h}$, as a combination of an encoding and a decoding procedure.  The lower dimensional
representation is obtained through an \textit{encoder} architecture
$\mathcal{F}_{\text{en}}:\boldsymbol{u}_{h}\to \boldsymbol{z}$, whose output is
subsequently mapped back to the high-dimensional space by a \textit{decoder}
network $\mathcal{F}_{\text{de}}:\boldsymbol{z} \to
\tilde{\boldsymbol{u}}_{h}$. The identity map is achieved by the composition of
the encoder and decoder, $\mathcal{I} = \mathcal{F}_{\text{de}}\circ
\mathcal{F}_{\text{en}}$,  allowing the autoencoder to simultaneously learn a compressed
representation of the input and its reconstruction.

Through a training process, the autoencoder's parameters are optimized to
minimize a suitable loss function, which is designed to measure the discrepancy
between the reconstructed output $\tilde{\boldsymbol{u}}_{h}$ and the original
input $\boldsymbol{u}_{h}$.

In our approach, we leverage the Sinkhorn divergence, a proxy for the
Wasserstein distance, as the chosen loss function. This decision offers
numerous advantages compared to more traditional options like mean squared
error or cross-entropy. These benefits include improved stability during
training, resilience against overfitting, robustness to noise, and faster
convergence \cite{frogner2015learning}.

Furthermore, we introduce an additional penalty term to enforce the encoded
lower-dimensional representation to coincide with the kPOD
reduced one. This ensures that the forward map is approximated as $F\approx
\mathcal{F}_{\text{en}}$, and the backward map as $F^{-1}\approx
\mathcal{F}_{\text{de}}$.  This approach provides a consistent nonlinear
generalization of POD in two aspects: the compression-decompression through the
autoencoder and the nonlinear projection on the kPOD principal components.

The main limitation of the use of autoencoders is the high number of parameters
to train in the case of systems with high-dimensional states. In fact, the
number of parameters of fully connected autoencoders exceeds by far the size of
the state itself. For this reason, we have also investigated the use of
\textit{Convolutional Autoencoders} (CAE) \cite{cdnn}, which have layers
characterized by the sharing of parameters and local connectivity, allowing for
a considerable reduction of the trainable parameters w.r.t.\ feedforward
autoencoders.

Overall, our study introduces a new framework for model order reduction of
partial differential equations, combining KPOD for nonlinear reduction, DL for
the inverse map approximation, and OT for kernel function definition and neural networks training.
\section{Related Works}
\label{sec:related}

Reduced order modeling techniques have been extensively studied in the scientific computing community due to their ability to efficiently reduce the computational cost of solving complex problems modeled by PDEs. In particular, machine learning-based ROMs \cite{brunton2019} have gained significant attention in recent years owing to their flexibility in handling high-dimensional and nonlinear systems.
Several deep learning architectures, including autoencoders \cite{autoencoder}, convolutional neural networks (CNNs) \cite{cdnn}, and recurrent neural networks (RNNs) \cite{rnn}, have been applied to model order reduction for scientific applications.

In the works by Fresca et al. \cite{fresca2021,fresca2022}, a non-intrusive approach to overcome common limitations of conventional ROMs is proposed by using deep learning-based reduced order models (DL-ROMs) that learn both the nonlinear trial manifold and the reduced dynamics. They also propose a possible way to avoid an expensive training stage of DL-ROMs by performing a prior dimensionality reduction through POD and relying on a multi-fidelity pretraining stage. A further development has been introduced in \cite{PichiGraphConvolutionalAutoencoder2023} by Pichi et al.\ exploiting a graph convolutional autoencoder (GCA-ROM) to consistently treat complex and possibly varying domains, and reduce the amount of training data needed thanks to the geometric bias, resulting in a light and interpretable architecture.

Although DL-ROMs have shown great promise in reducing computational costs, other researchers have pursued alternative methods for addressing the limitations of linear-subspace ROMs. One such approach is the projection of dynamical systems onto nonlinear manifolds, as proposed by Lee et al.\ \cite{Lee2020}. The authors used minimum-residual formulations at the time-continuous and time-discrete levels, and they could significantly outperform even the optimal linear-subspace ROM on advection-dominated problems. Building on this work, Romor et al.\ \cite{romor2023} have developed a nonlinear manifold method that achieves hyper-reduction through reduced over-collocation and teacher-student training of a reduced decoder. This approach offers a promising alternative to purely data-driven methods. Additionally, M\"{u}ckeet al.\ \cite{Mcke2021} present a non-intrusive ROM approach based on model learning, consisting of a nonlinear dimensionality reduction stage using convolutional autoencoders and a parameterized time-stepping stage using memory-aware neural networks.

While many researchers have focused on developing ROMs that reduce the computational cost of existing high-fidelity simulations, others have explored the possibility of using DL frameworks to completely replace the high-fidelity model, as discussed in  Karniadakis et al.\ \cite{karniadakis2021}: in Raissi et al.\ \cite{raissi2019,raissi2018}, physics-informed neural networks (PINNs) are introduced, which are trained to solve supervised learning tasks while incorporating given laws of physics described by PDEs.

Furthermore, there has been a growing interest in combining classical ROM techniques \cite{benner2017model, benner2017, Grepl2007,schilders,vol1,vol2}, such as POD \cite{pod}, with machine learning methods \cite{machine-learning-book}. In particular, Ubbiali et al.\ \cite{ubbiali} proposed the POD-NN method, which utilizes a neural network for regressing the reduced order coefficient instead of projecting the governing equations on the POD modes. This technique has since then been employed for investigations in various physical contexts \cite{hesthaven2019,aerostructural-podnn,wang2019,Pichi2021}.

Despite the promising results of POD-NN in reducing the computational cost of reduced order modeling, there are some drawbacks, such as overfitting and the requirement for a significant amount of data. As a result, other researchers have explored the use of kernel-based methods, such as kernel POD (kPOD) \cite{kpca} and kernel CCA (KCCA) \cite{kcca}, which offer an alternative approach with less reliance on the availability of large training sets.
In the study by Díez et al. \cite{Dez2021}, a proposed methodology for reduced order modeling suggests using kPOD to discover a lower-dimensional nonlinear manifold suitable for solving parametric problems. The kernel employed in kPOD is defined based on physical considerations, and the approximation space is enriched with cross-products of the snapshots.
An alternative approach proposed by Salvador et al. \cite{salvador} leverages an NN architecture to find the reduced basis coefficients, which are subsequently used to recover a linear backward mapping.

Kernel-based methods have been successful in reducing the dimensionality of high-dimensional data, but they may fall short in capturing the geometric features of the underlying data distribution. Optimal transport (OT) theory \cite{santambrogio2015optimal}, on the other hand, focuses on the geometry of the data distribution and has become an increasingly popular tool in machine learning \cite{peyre2019computational}. One of its main applications is the use of the Wasserstein distance to measure dissimilarity between high-dimensional data \cite{frogner2015learning}. The Sinkhorn algorithm, popularized in the OT context by Cuturi \cite{cuturi13}, efficiently computes the regularised Wasserstein distance and has been used in numerous works, including the Wasserstein GAN \cite{wasserstein-gan}. This method uses the Wasserstein distance to measure the discrepancy between true and generated data in a generative adversarial network.

In the context of ROMs for parametrized PDEs, Mula et al. \cite{mula2020} addressed the reduction problem using a nonlinear approximation in the general metric setting given by the $L^2$-Wasserstein space. Theoretical and numerical results are showcased for one-dimensional conservative PDEs, and the approach is extended to non-conservative problems in \cite{mula2022}. In \cite{mula2023}, the authors propose a sparse approximation and structured prediction approach using Wasserstein barycenters introduced by Agueh et. al. \cite{agueh2011barycenters}. The method is adaptive, mass preserving, and has the potential to overcome known bottlenecks of classical linear methods for hyperbolic conservation laws.

In the task of parametric interpolation, Torregorsa et al. \cite{Torregrosa2022} presented a non-intrusive parametric OT-based model using the sparse Proper Generalized Decomposition (sPGD) regression. On the other hand, Taddei et al. \cite{taddei2022} proposed an interpolation technique for parametric advection-dominated problems based on OT of Gaussian models. This promising approach identifies the coherent structure to be tracked using a maximum likelihood estimator, relying on nonlinear interpolation through OT maps between Gaussian distributions. In addition, in \cite{blickhan2023registration}, the author proposes a registration method \cite{Taddei2020} that aligns solution features to facilitate approximation using classical reduced basis methods, employing OT mappings in a purely data-driven manner.

Overall, the combination of deep learning and optimal transport theory, together with the use of custom kernels, has the potential to yield highly effective and efficient ROM techniques, which can pave the way for many possible methodological developments and complex applications across various scientific domains.
\section{Optimal Transport Theory}
\label{sec:ot}
Optimal Transport (OT) theory \cite{santambrogio2015optimal} is a rigorous mathematical framework used to solve transportation problems. It finds applications in diverse fields such as probability theory, economics, image processing, and machine learning. In our framework, OT theory plays a crucial role by leveraging the metric it establishes, known as the Wasserstein distance. This metric serves two primary purposes in our approach:

\begin{itemize}
    \item Dimensionality reduction via OT kernel: the Wasserstein distance is employed in combination with kPOD to perform non-linear projections, facilitating the parametric regression task and enhancing the modeling capabilities of our framework.
    \item Training neural networks: we incorporate the Wasserstein metric  into the learning process to approximate the kPOD backward map, allowing the neural networks to easily capture the underlying features and to enhance its accuracy.
\end{itemize}

In the following, we present a brief introduction  OT theory, providing the fundamental principles underlying the computation of the Wasserstein distance exploited within our architecture. The focus will be on the discrete case, which is of particular importance for our study.

\subsection{Formulation}

Originated as the problem of finding the optimal way of moving a pile of soil between two locations while minimizing the total transportation cost, Monge formulated the OT theory seeking the optimal mapping between two probability measures, denoted as $\mu$ and $\nu$. Let $\mathcal{X}$ and $\mathcal{Y}$ be the source and target metric spaces respectively, and $c:\mathcal{X}\times\mathcal{Y}\rightarrow \mathbb{R}^+$ be the cost function that quantifies the transportation cost from $\mathcal{X}$ to $\mathcal{Y}$.

In addition, given a measurable function $T: \mathcal{X}\rightarrow \mathcal{Y}$, the push-forward of $\mu$ by $T$ is defined as the measure $T_\#\mu$ on $\mathcal{Y}$ such that $T_\#\mu(A) = \mu(T^{-1}(A))$ for any measurable set $A \subset \mathcal{Y}$. We refer to $T$ as a \textit{transport map} if $T_\#\mu = \nu$, \ie, if the transport map enables to express measures on the target space in terms of their pre-images.
The Monge problem can be formulated as finding the optimal transport map $T^*$ that minimizes the total transportation cost:

\begin{equation}
    T^* = \argmin_{T \ \text{s.t.}\ T_\#\mu = \nu} \int_{\mathcal{X}} c(x ; T(x)) d \mu .
\end{equation}
Despite its mathematical clarity, this formulation presents several challenges: (i) the functional to be minimized is completely non-linear and asymmetric with respect to $\mathcal{X}$ and $\mathcal{Y}$, and (ii) the problem is ill-posed, in the sense that the set of admissible transport maps could be empty.

To overcome the limitations of the Monge formulation, Kantorovich introduced an alternative approach by considering couplings between measures instead of direct mapping. A probability measure $\pi \in \mathcal P(\mathcal X \times \mathcal Y)$ is called a \textit{transport plan} if $(p_\mathcal{X})_\#\pi = \mu$ and $(p_\mathcal{Y})_\#\pi = \nu$. Here, $p_\mathcal{X}: \mathcal{X} \times \mathcal{Y} \rightarrow \mathcal{X}$ and $p_\mathcal{Y}: \mathcal{X} \times \mathcal{Y} \rightarrow \mathcal{Y}$ denote the coordinate projections that extract the first and second components of $(x,y) \in \mathcal{X}\times\mathcal{Y}$, respectively.

We will indicate with $\Pi(\mu,\nu)$  the set of all admissible transport plans between $\mu$ and $\nu$. We can notice how a transport plan $\pi$ is a probability measure on the product space $\mathcal{X}\times \mathcal{Y}$ with $\mu$ and $\nu$ as its marginals, that specify the amount of mass to be transported from each point in $\mathcal{X}$ to each point in $\mathcal{Y}$.

The Kantorovich problem consists in finding the optimal transport plan $\pi^*\in\Pi(\mu,\nu)$ that minimizes the total cost:

\begin{equation}
    \pi ^* = \argmin_{\pi \in \Pi(\mu,\nu)}\int_{\mathcal{X}\times \mathcal{Y}}c(x, y)d\pi(x,y).
\end{equation}
A key advantage of such formulation is that it allows for the existence of partial transportations where mass can be split or merged during the transportation process, \ie, without the requirement of a bijection between the measures, which was implicit in the definition of $T_\# \mu$. To quantify the optimal transportation cost between two measures $\mu$ and $\nu$, we introduce the \textit{Wasserstein distance}:
\begin{equation}
    W_p(\mu, \nu) = \left( \min_{\pi \in \Pi (\mu,\nu)} \int_{\mathcal{X} \times \mathcal{Y}} c(x,y)^p d\pi(x,y) \right)^{\frac{1}{p}},
\end{equation}
where $p \geq 1$ is a parameter that determines the order of the distance.  The Wasserstein distance provides a quantitative measure of the optimal transportation cost between two measures, taking into account the underlying cost function and the
geometry of the measures.

\begin{figure}[h]
    \centering
    \includegraphics[width=0.5\linewidth]{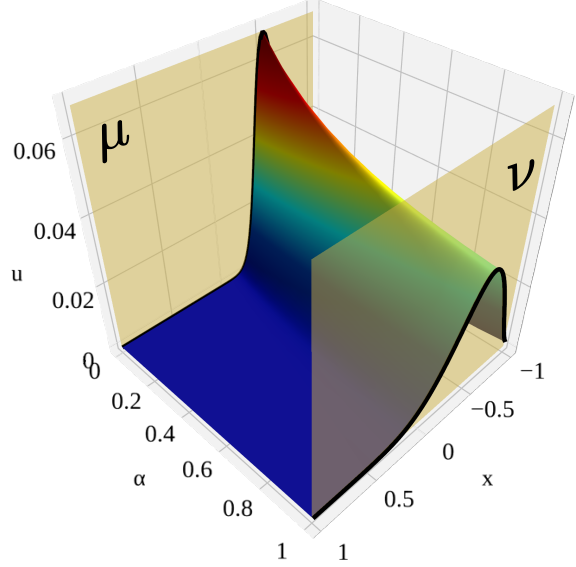}
    \caption{Wasserstein barycenters between the solution to the Burgers equation at two different time instants. The parameter $\alpha$ represents the interpolation parameter between the two fields $\mu$ and $\nu$, which are normalized to treat them as probability distributions.}
    \label{fig:barycenter}
\end{figure}

In addition, the Wasserstein distance allows for the definition of the Wasserstein barycenter \cite{agueh2011barycenters}, which enables obtaining a representative measure that captures the central tendency or average distribution from a given set of probability measures (see Figure \ref{fig:barycenter}).
In the context of reduced-order modeling for parametric PDEs, this tool can be employed to construct interpolation strategies, enabling structured prediction and metamodelling for parametrized families of measures \cite{mula2023}.

\subsection{Discrete Optimal Transport}

In the discrete case, the Kantorovich problem can be formulated as a linear programming problem, which can be efficiently solved using various optimization techniques, such as the simplex algorithm \cite{transportation-simplex}. In this case, $\mu$ and $\nu$ are represented as discrete probability measures with finite supports:

\begin{equation} \begin{aligned} & \mu=\sum_{i=1}^n \boldsymbol{a}_i
    \delta_{x_i} , \quad \boldsymbol{a} \in \mathbb{R}_+^{n} \quad \text{s.t.}
    \quad \boldsymbol{a}^{T} \mathbf{1}_n=1, \\ & \nu=\sum_{j=1}^m
    \boldsymbol{b}_j \delta_{y_j},\quad \boldsymbol{b} \in \mathbb{R}_+^{m}
    \quad \text{s.t.} \quad \boldsymbol{b}^{T} \mathbf{1}_m =1, \end{aligned}
\end{equation}
where $\delta_{x_i}$ and $\delta_{y_j}$ denote Dirac delta functions at points $x_i$ and $y_j$, respectively, and $\boldsymbol{a}$ and $\boldsymbol{b}$ are non-negative vectors representing the weights of these points.  The notation $\boldsymbol{1}_{n}= [1, \ldots, 1] \in \mathbb{R}_{+}^{n}$ refers to the uniform histogram.

The cost function $c$ is encoded in the cost matrix $\mathbf{C} \in \mathbb{R}_{+}^{n \times m}$, where $n$ and $m$ are the sizes of the supports of $\mu$ and $\nu$, respectively. Each entry  $(\mathbf{C})_{ij}$ represents the cost of transporting mass from the $i$th point in $\mathcal{X}$ to the $j$th point in $\mathcal{Y}$, raised to the $p$ power:
\begin{equation}
    (\mathbf{C})_{ij} =c\left(x_i, {y}_j\right)^p.
\end{equation}
A generic transport plan is then given by a coupling\footnote{In the discrete case, the probabilities are encoded by vectors or matrices for couplings; this distinction is here highlighted through the boldface notation.} $\boldsymbol{\pi} \in \mathbb{R}_+^{n\times m}$, which belongs to the following set:

\begin{equation} \Pi_{\text{dis}}(\mu, \nu) {=}\left\{\mathbf{P} \in
    \mathbb{R}_{+}^{n \times m} \mid \mathbf{P} \mathbf{1}_m=\boldsymbol{a}, \
\mathbf{P}^T \mathbf{1}_n=\boldsymbol{b}\right\}.  \end{equation}
Clearly $\boldsymbol{\pi} \in \Pi_{\text{dis}}(\mu,\nu)$ is a constraint ensuring that
$\boldsymbol{\pi}$ is a valid transport plan with $\mu$ and $\nu$ as marginals,
in fact: \begin{equation} \boldsymbol{\pi} \boldsymbol{1_n} = \boldsymbol{\pi}
    \left[\begin{array}{c}
1 \\
1 \\
\vdots \\
1 \\
\end{array}\right]=\left[\begin{array}{cc}
\sum_j \pi_{1j}  \\ \sum_j \pi_{2j}  \\
\vdots \\
\sum_j \pi_{nj} \\
\end{array}\right] = \boldsymbol{a}, \quad
\boldsymbol{\pi}^{T} \boldsymbol{1_m} =
\boldsymbol{\pi}^{T} \left[\begin{array}{c}
1 \\
1 \\
\vdots \\
1 \\
\end{array}\right]=\left[\begin{array}{cc}
\sum_j \pi_{j1}  \\ \sum_j \pi_{j2}  \\
\vdots \\
\sum_j \pi_{jm} \\
\end{array}\right] = \boldsymbol{b}.
\end{equation}
The discrete counterpart of the Wasserstein distance is thus given by:
\begin{equation}
    W_p^p({\mu},{\nu})=\min _{\substack{\mathbf{P} \in \Pi_{\text{dis}}(\mu,\nu)}}\left\langle\mathbf{P}, \mathbf{C}\right\rangle,
\end{equation}
where $\left\langle \cdot, \cdot\right\rangle$ represents the Frobenius inner product between matrices, defined as the sum of the products of the corresponding components.

The advantage of the discrete Kantorovich problem is that it can be efficiently solved using linear programming techniques. The transportation simplex algorithm \cite{transportation-simplex}, for example, has an average polynomial complexity in the number of variables and constraints but exponential worst-case complexity. Other methods, such as the auction algorithm \cite{auction}, achieve a worst-case cubic complexity, which may be prohibitively costly for large-scale problems, especially when dealing with high-dimensional data or large datasets \cite{ot-cost}. Various optimization techniques have been proposed to speed up the computation, such as the entropic regularization method, which will be discussed in the next section.

\subsection{Entropic Regularization and Sinkhorn Algorithm}
\label{subsec:sink}

Computing the Wasserstein distance is a nontrivial task, especially for high-dimensional data. However, the Sinkhorn algorithm, exploiting entropic regularization, provides an efficient and scalable method for its computation.

\begin{figure}
    \centering
    \includegraphics[width=1\textwidth]{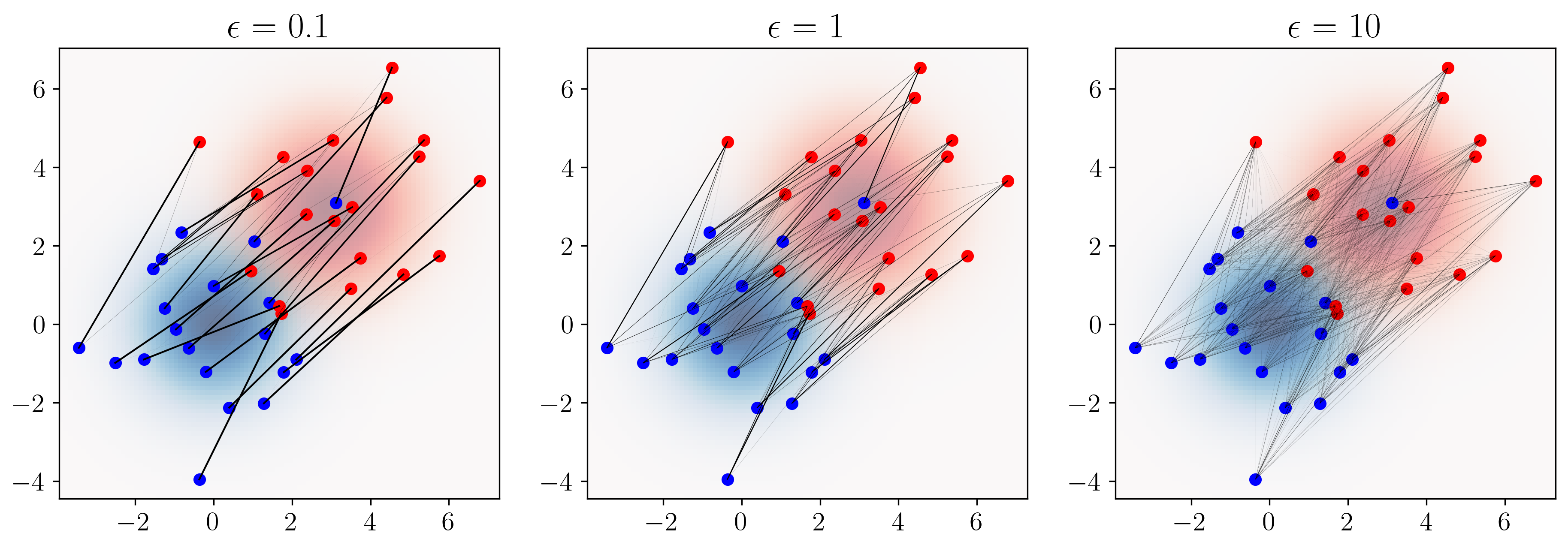}
    \caption{Optimal transport plans between two Gaussian distributions with
    different regularization parameters. The blue and red dots represent
samples drawn from the two distributions, and the black lines indicate the
transport plans between them, with thicker lines indicating higher plan values.
As the regularization parameter increases, the transport plans become smoother
and more diffuse. }%
    \label{fig:sink_}
\end{figure}

The entropic regularization introduces a small amount of entropy into the optimization problem, leading to a differentiable and convex optimization problem that can be efficiently solved using iterative methods. This procedure involves the introduction of a parameter $\epsilon > 0$ that controls the amount of regularization (see Figure \ref{fig:sink_}).

The regularized Wasserstein distance is given by:
\begin{equation}
    W_{\text{reg}}(\mu, \nu) \stackrel{\text { def }}{=} \min_{\boldsymbol{\pi} \in \Pi_{\text{dis}}(\mu, \nu)} \langle \mathbf{C}, \boldsymbol{\pi} \rangle - \epsilon H(\boldsymbol{\pi}),
\label{eq:wasserstein_regularized}
\end{equation}
where $H(\pi)$ is the \textit{Shannon entropy} of the transport plan $\pi$, defined as:

\begin{equation}
    H(\boldsymbol{\pi}) \stackrel{\text { def }}{=} - \sum_{i = 1} ^n \sum _{j = 1}^m \boldsymbol{\pi}_{ij}
 (\log(\boldsymbol{\pi}_{ij})-1).
 \label{eq:entropy}
\end{equation}
The entropic regularization introduces the smoothing term
$H(\boldsymbol{\pi})$, which is uniformly concave with respect to its argument.

This allows for more efficient optimization since the regularized Problem
\eqref{eq:wasserstein_regularized} is uniformly convex and thus admits a unique
minimum.
To achieve the optimum for Problem \eqref{eq:wasserstein_regularized}, a first simplification lies in writing the corresponding Karush-Kuhn-Tucker \cite{Rockafellar1970}
conditions. Doing so yields
the following characterization:
\begin{equation}
 \boldsymbol{\pi}^* \text { solution to Problem \eqref{eq:wasserstein_regularized} } \Leftrightarrow\left\{\begin{array}{l}
        \boldsymbol{\pi}^* \in \Pi_{\text{dis}}(\mu,\nu) , \\
  \mathbf{K} \in \mathbb{R}_+^{n \times m}: \quad \mathbf{K}_{ij} = e^{\frac{\mathbf{C}_{ij}}{\epsilon}},\\
        \exists\ \boldsymbol{u} \in \mathbb{R}^n,\ \boldsymbol{v} \in \mathbb{R}^m: \quad \boldsymbol{\pi}^{*}_{ij}=\boldsymbol{u}_i \mathbf{K}_{ij} \boldsymbol{v}_j.
\end{array} \right.
\label{ktt-conditions}
\end{equation}
Thanks to the equivalence provided by Equation \eqref{ktt-conditions}, one can
solve the Kantorovich problems using the Sinkhorn-Knopp's fixed point iteration
\cite{sink} reported in Algorithm \eqref{alg-sink}.  The algorithm works by
iteratively scaling the rows and columns of the transport plan until
convergence. The resulting sequence of transport plans is guaranteed to
converge to the optimal one.

\begin{algorithm}[h]
\SetAlgoLined
  \SetKwInput{Input}{Input}
  \SetKwInput{Output}{Output}
  \Input{Two probability vectors $\boldsymbol{a} \in \mathbb{R}^n$, $\boldsymbol{b} \in \mathbb{R}^m$, a cost matrix $C \in \mathbb{R}^{n \times m}$, and a regularization parameter $\epsilon > 0$}
  \Output{Optimal transport plan $\mathbf{P}$ of size $n \times m$}
  Initialize $\boldsymbol{u}^{(0)} \leftarrow \boldsymbol{1}_n$ and $\boldsymbol{v}^{(0)} \leftarrow \boldsymbol{1}_m$\;
  Set $\ell \leftarrow 0$\;
  \tcp{Compute the matrix $\mathbf{K}$}
  $(\mathbf{K})_{ij} \leftarrow e^{-\frac{(\mathbf{C})_{ij}}{\epsilon}}$ for all $i$ and $j$\;
  \While{not converged}{
    \tcp{Update dual variables $u$ and $v$}
    $\boldsymbol{u}^{(\ell+1)} \leftarrow \frac{\boldsymbol{a}}{\mathbf{K} \boldsymbol{v}^{(\ell)}}$\;
    $\boldsymbol{v}^{(\ell+1)} \leftarrow \frac{\boldsymbol{b}}{\mathbf{K}^{\top} \boldsymbol{u}^{(\ell+1)}}$\;
    $\ell \leftarrow \ell + 1$\;
  }
  \tcp{Compute the regularized transport plan $P$ }
  $(\mathbf{P})_{ij} \leftarrow (\boldsymbol{a}^{(\ell)})_i (\mathbf{K})_{ij} (\boldsymbol{b}^{(\ell)})_j$ for all $i$ and $j$\;
  \Return $\mathbf{P}$\;
\mycaption{Sinkhorn-Knopp's Fixed Point Iteration for Optimal Transport}
\label{alg-sink}
\end{algorithm}

Algorithm \eqref{alg-sink} is iterated until convergence is achieved, as
measured by the duality gap or another stopping criterion
\cite{peyre2019computational}. The output is the regularized transport plan
$\boldsymbol{\pi}_\epsilon$ which can be used to compute the \textit{Sinkhorn
distance} \cite{cuturi13}:
\begin{equation}
W_{\epsilon}(\mu, \nu)\stackrel{\text { def }}{=}\ \langle \boldsymbol{\pi}_\epsilon,  \mathbf{C}  \rangle.
\label{eq:sinkhorn-distance}
\end{equation}
The issue with the former quantity is that, unlike the Wasserstein distance, it
can be proven to not be positive and definite. In particular, for any given target
measure $\nu$, there exist a measure $\mu$ s.t.:
\begin{equation}
W_{\epsilon}(\mu, \nu)< W_{\epsilon}(\nu, \nu)\neq 0 .
\end{equation}
However, we can introduce additional corrective terms to debias $W_{\epsilon}$, obtaining the so-called \textit{Sinkhorn divergence} \cite{Ramdas2017,genevay2018}:
\begin{equation}
    S_{\epsilon}(\mu, \nu)\stackrel{\text { def }}{=}\ W_{\epsilon}(\nu, \mu) - \frac{1}{2}W_{\epsilon}(\nu, \mu) - \frac{1}{2}W_{\epsilon}(\nu, \mu).
\end{equation}
The previous quantity trivially satisfies $S_{\epsilon}(\mu,\mu ) = 0$, and provides a valid dissimilarity measure between two probability distributions \cite{feydy2019interpolating}.
Furthermore, since the computation of $S_{\epsilon}(\mu, \nu)$  only involves matrix-vector multiplications, it is highly parallelizable and can take advantage of GPUs, resulting in a significant reduction in computational cost. Algorithm \eqref{alg-sink} can be further optimized to achieve a computational complexity of $\mathcal{O}(n \log n + m \log m)$ with a linear memory footprint \cite{feydy-thesis}, making it suitable for solving large-scale optimization problems.
\section{Kernel Proper Orthogonal Decomposition (kPOD)}
\label{sec:kpod}

The Kernel Proper Orthogonal Decomposition (kPOD) algorithm \cite{kpca} is a kernel-based extension of the classical POD algorithm.
kPOD leverages the idea that datasets with nonlinearly separable attributes can be embedded into a higher-dimensional feature space, where the information becomes more spread out and easier to compress through POD.
When projected back to the original space, the resulting modes yield a nonlinear projection that better represents the data. To avoid the explicit calculation of the higher-dimensional coordinates, kPOD uses the kernel function to implicitly map the data to the feature space, reducing the method's computational complexity. As a result, kPOD can perform the computation in the original space while still benefiting from the advantages of the nonlinear projection.

We will describe the method starting from the description of POD, which relies on the Singular Value Decomposition (SVD) of the snapshot matrix $\mathbf{S}$ defined in \eqref{eq:snap}. The SVD provides a factorization of $\mathbf{S} \in \mathbb{R}^{N_{h}\times N_{s}}$ into the form:
\begin{equation}
 \mathbf{S} = \mathbf{U} \mathbf{\Sigma} \mathbf{V^{T}},
\end{equation}
 where $\mathbf{U} \in \mathbb{R}^{N_{h}\times N_{h}}$ and $\mathbf{V} \in \mathbb{R}^{N_{s}\times N_{s}}$ are orthogonal matrices, and $\mathbf{\Sigma} \in \mathbb{R}^{N_{h}\times N_{s}}$ is a diagonal matrix.

 Let $\mathbf{U}^\star = [\boldsymbol{\varphi}_1 | \boldsymbol{\varphi}_2| \ \cdots| \ \boldsymbol{\varphi}_k] \in \mathbb{R}^{N_h \times k}$ be the truncated matrix consisting of the first $k$ columns of $\mathbf{U}$. The low-dimensional representation of a snapshot is given by $\boldsymbol{z} = \mathbf{U^{\star T}} \boldsymbol{u}_{h} \in \mathbb{R}^k$, so that the samples in $\mathbf{S}$ are mapped to $\mathbf{Z} = \mathbf{U^{\star T}} \mathbf{S} \in \mathbb{R}^{k \times N_s}$. The reduced order model is interpreted as looking for an approximation to $\boldsymbol{u}_{h}$ in the $k$-dimensional linear subspace generated by the basis $\{\boldsymbol{\varphi}_j\}_{j = 1}^{k}$. Notice that the SVD of the snapshot matrix $\mathbf{S}$ directly provides a diagonalization of the Gram matrix $\mathbf{G} = \mathbf{S^T} \mathbf{S} \in \mathbb{R}^{N_s \times N_s}$ as follows:

\begin{equation}
    \label{eig-prob-G}
    \mathbf{G} = \mathbf{V}\mathbf{[\Sigma ^T \Sigma]V ^ T = V \Lambda V ^T},
\end{equation}
where $\mathbf{\Lambda}$ is a diagonal matrix with the eigenvalues of $\mathbf{G}$ on its diagonal. In particular, the compact expression for constructing $\mathbf{Z}$  is given in terms of the Gram matrix as follows:

\begin{equation}
\label{kpod-reduction}
\mathbf{Z=V ^{\star  T} G} ,
\end{equation}
where $\mathbf{V}^\star \in \mathbb{R}^{N_s \times k}$ is the matrix consisting of the first $k$ columns of $\mathbf{V}$.

To apply kPOD, we first introduce an arbitrary transformation $\Phi: \mathbb{R}^{N_{h}} \rightarrow \mathbb{R}^{N_{H}}$ that maps the original data to a higher-dimensional feature space, where nonlinearities may become linear. Let us assume to have $N_s$ transformed snapshots $\{\Phi(\boldsymbol{u}^{j}_{h})\}_{j = 1}^{N_{s}}$, from which we can define the matrix $\tilde{\mathbf{S}} = [\Phi(\boldsymbol{u}^1_{h})| \Phi(\boldsymbol{u}^2_{h})| ... | \Phi(\boldsymbol{u}^{N_s}_{h})] \in \mathbb{R}^{N_H \times N_s}$. We aim at performing POD on $\tilde{\mathbf{S}}$ to extract a low-dimensional linear subspace that captures most of the information in the transformed snapshots.

However, computing the SVD of $\tilde{\mathbf{S}}$ directly may be computationally prohibitive if $N_{H}$ is very large, and we do not know the optimal choice of $\Phi$ that would allow us to untangle the underlying nonlinear manifold. Instead, we can compute the transformed Gram matrix $\tilde{\mathbf{G}}$ by introducing a bivariate symmetric form $\kappa(\cdot,\cdot): \mathbb{R}^{N_{h}} \times \mathbb{R}^{N_{h}} \rightarrow \mathbb{R}^k$ called \textit{kernel function} that acts as a proxy for the inner product in the high-dimensional feature space, but takes as argument the pre-image of the transformed samples:

\begin{equation}
    [\tilde{\mathbf{G}}]_{ij} = \langle \Phi(\boldsymbol{u}^{i}_{h}), \Phi(\boldsymbol{u}^{j}_{h}) \rangle = \kappa(\boldsymbol{u}^{i}_{h},\boldsymbol{u}^{j}_{h}).
\end{equation}
By doing so, we can still rely on \eqref{eig-prob-G} to obtain $\tilde{\mathbf{V}}^{\star}$ and use \eqref{kpod-reduction} to perform the model order reduction task.
In particular, we can introduce the so-called \textit{forward mapping} F:

\begin{equation}
\label{forward-map}
\begin{aligned}
F: \mathbb{R}^{N_{h}} & \rightarrow \mathbb{R}^k, \\
{\boldsymbol{u}_{h}} & \mapsto {\boldsymbol{z}} = {\tilde{\mathbf{V}}}^{\mathbf{\star T}}  {g}({\boldsymbol{u}_{h}}), \qquad {g}({\boldsymbol{u}_{h}}) = [\kappa({\boldsymbol{u}_{h}}, {\boldsymbol{u}^{i}_{h}})]_{i=1}^{N_s} \ .
\end{aligned}
\end{equation}
However, we do not have an explicit expression for the kPOD modes, which would be contained in the unknown matrix $\tilde{\mathbf{U}}$. This represents a significant problem because having only $\tilde{\mathbf{V}}^{\star}$ does not allow us to correctly map backward an element $\boldsymbol{z} \in \mathbb{R}^k$ into its pre-image $\tilde{\boldsymbol{u}}_{h} \in \mathbb{R}^{N_{h}}$.
For a more detailed description of kPOD and the various strategies for defining a backward map, we recommend interested readers to refer to \cite{kpca,Dez2021}.

A crucial aspect consists in determining the backward mapping of the low-dimensional manifold with kPOD, which is less trivial than in the standard POD framework. In particular, defining the inverse of the forward mapping $F$ requires the selection of a proper target space $\mathcal{V} \subset \mathbb{R}^{N_{h}}$ where pre-images are sought. The \textit{backward mapping} is then defined as:

\begin{equation}
    \tilde{\boldsymbol{u}}_{h}=F^{-1}(\boldsymbol{z})\stackrel{\text { def }}{=}\arg \min _{\boldsymbol{u}_{h} \in \mathcal{V}}\|F(\boldsymbol{u}_{h})-\boldsymbol{z}\|.
\end{equation}
The definition of the backward mapping depends on the selection of the manifold $\mathcal{V}$ and the approximation criterion, which is the choice of the norm $\|\cdot\|$. We remark that, since $F$ is not injective, different $\tilde{\boldsymbol{u}}_{h} \in \mathbb{R}^{N_{h}}$ can have the same image $\boldsymbol{z}$, thus, the target space must be a proper low-dimensional manifold for the backward map to be well-defined.

To address this problem, we propose to use a regression approach based on Neural Networks (NN) \cite{goodfellow,freeman1991neural}. Specifically, we will use autoencoders \cite{autoencoder}, which have been shown to be effective in solving the inverse problem originated from kPOD.
\section{OT-based DL-framework}\label{sec:complete_framework}
In this section, we exploit the methodologies discussed in the previous sections to introduce the optimal transport based model order reduction approach we developed, which takes advantage of the Sinkhorn divergence kernel for the kPOD and of a neural network regression for the backward map problem.

\subsection{Sinkhorn Kernel}
Let $\mathcal{X} = \{x_{i}\}_{i = 1}^{N_{h}}$ be a set of discrete input data, and $\mu, \nu$ be two probability measures over $\mathcal{X}$, represented by the two probability mass functions $f,g:\mathcal{X} \to [0,1]$.

The kernel we propose for the kPOD is given by:
\begin{equation}
    \label{eq:sink-kern}
\kappa( \mu, \nu) = \left[\frac{1}{2} \sum_{i=1}^{N_{h}} {x_{i}}^{T} {x_{i}} f({x_{i}}) +
\frac{1}{2} \sum_{j=1}^{N_{h}} {y_{j}}^{T} {y_{j}} g({y_{j}}) \right]
-S_{\epsilon}(\mu,\nu),
\end{equation}
where $x_{i}$ and $y_{j}$ denote the support points of $\mu$ and $\nu$ respectively, and $S_{\epsilon}(\mu,\nu)$ can be efficiently computed using the Sinkhorn-Knopp algorithm (see Section \ref{sec:ot}).
The motivation behind this formula is to obtain a weighted version of the $l^{2}$ discrete product that takes into account the transport plan.

Considering the regularity of Monge's Problem, we ensure that an optimal transport map $T:\mathcal{X} \to \mathcal{X}$ exists, and it is bijective. Consequently, the inverse map $T^{-1}:\mathcal{X} \to \mathcal{X}$ also exists, leading to the relationships ${x}_{i} = T^{-1}({y}_{i})$ and ${y}_{i} = T({x}_{i})$. Here, the set $\left\{y_{i}\right\}_{i = 1}^{N_{h}}$ represents a proper permutation of $\left\{x_{i}\right\}_{i = 1}^{N_{h}}$ which solves the Monge's Problem.

Exploiting these properties, one can easily rewrite the bivariate symmetric form as:

\begin{equation}
    \label{kernel_form2}
\kappa(\mu, \nu) = \sum_{i=1}^{N_{h}} {x_{i}}^{T} T({x_{i}}) f({x_{i}}) =
\sum_{i=1}^{N_{h}} {y_{i}}^{T} T^{-1}({y_{i}}) g({y_{i}}) = \kappa (\nu, \mu) .
\end{equation}
To prove the equivalence between Equations \eqref{eq:sink-kern} and \eqref{kernel_form2}, we start by considering the case of computing $\kappa(\mu,\mu)$. In this scenario, the optimal transport map is equal to the identity, denoted as $\mathcal{I}:\mathcal{X} \to \mathcal{X}$. Consequently, we have:

\begin{equation}
\kappa(\mu, \mu) = \sum_{i = 1}^{N_{h}} {x_{i}}^{T} \mathcal{I}({x_{i}}) f({x_{i}}) =
\sum_{i = 1}^{N_{h}} {x_{i}}^{T} {x_{i}} f({x_{i}}).
\end{equation}
This reduces the kernel $\kappa(\mu, \mu)$ to the second moment of the probability distribution $\mu$.

Next, we develop the Wasserstein distance in the quadratic cost case:
\begin{equation}
    W_{2}^{2}(\mu,\nu) = \sum _{ i = 1 }^{N_{h}}   \mid{x_{i}} - T({x_{i}})\mid^{2} f({x_{i}}) = \sum _{i =1 }^{N_{h}}( {x}_{i}^{T} {x_{i}} + T({x}_{i})^{T} T({x}_{i}) - 2 {x_{i}}^{T}T({x}_{i})) f({x}_{i}).
\end{equation}
After isolating the mixed product term, which corresponds to the kernel $\kappa(\mu,\nu)$, we can use the pushback $T_{\#} \mu= \nu$ to obtain:
\begin{equation}
\kappa( \mu, \nu) = \left[\frac{1}{2} \sum_{i=1}^{N_{h}} {x_{i}}^{T} {x_{i}} f({x_{i}}) +
\frac{1}{2} \sum_{j=1}^{N_{h}} {y_{i}}^{T} {y_{i}} g({y_{i}})\right]
-W_{2}^{2}(\mu,\nu) .
\end{equation}
Our approach aims to utilize the Sinkhorn divergence as an approximation for the $W_2^{2}(\cdot,\cdot)$ distance, thereby defining the kernel presented in Equation \eqref{eq:sink-kern}. As shown, this kernel is equivalent to the one presented in Equation \eqref{kernel_form2}, but it eliminates the need for an explicit computation of the optimal transport map $T$.

\begin{remark}
Consider a possible extension of the Gaussian kernel, where we replace the Euclidean distance with the Sinkhorn divergence:

\begin{equation}
\kappa( \mu, \nu) = \exp\left(-\frac{1}{\sigma} (S_{\epsilon}(\mu,\nu))^{2}\right),
\end{equation}
and the parameter $\sigma$ controls the smoothness of the kernel.
It can be proved that this kernel, known as the Wasserstein exponential kernel \cite{Plaen2020}, lead to a Gram matrix representing a first-order approximation of the kernel we propose, up to a rank-one matrix.
\end{remark}

\subsection{Neural network architectures}
To construct the backward map, we employ a neural network regression approach to learn the mapping between the reduced-order coefficients and the full-order solution, which corresponds to the inverse mapping of kPOD. We explore two distinct neural network architectures for this purpose: a decoder-only architecture and a complete autoencoder architecture.

\subsubsection{Decoder-based approach}

In the decoder-only architecture, we construct a neural network that takes the kPOD coefficients $\boldsymbol{z} \in \mathbb{R}^{k}$ as input and generates the full-order solution $\boldsymbol{u}^{D}_{\text{NN}} \in \mathbb{R}^{N_{h}}$ as output.

The decoder function $\psi^{D}(\boldsymbol{z};{\Theta}^{D})$ is defined by a set of weights and biases ${\Theta^{D}}$, which are learned during the training process. Once trained, it can be used to evaluate the backward mapping as:

\begin{equation}
      \boldsymbol{u}^{D}_{\text{NN}} = f^{D}(\boldsymbol{z};{\hat{\Theta}}^{D}) \approx F^{-1}(\boldsymbol{z}) .
\end{equation}
To find the optimal parameter vector $\hat{\Theta}^{D}$, we employ a supervised learning paradigm that involves optimizing a properly selected loss function $\mathcal{J}^{D}$:

\begin{equation}
\label{nn-opt}
{\hat{\Theta}^{D}}=\underset{{\Theta}^{D}}{\operatorname{argmin}} \ \mathcal{J}^{D}({\Theta}^{D}),
\end{equation}
which measures the discrepancy between the true full-order solutions and the output of the neural network. The optimization process aims to minimize this discrepancy by adjusting the parameters ${\Theta}^{D}$, leading to an accurate decoder model that can effectively reconstruct the full-order solutions from the kPOD coefficients.

\subsubsection{Autoencoder-based approach}
\label{sec:ae-based-approach}

In the complete autoencoder architecture, we extend the decoder-based approach by adding an encoder component to the neural network.
The encoder function $\psi^{E}(\cdot;{\Theta^{E}})$ is parameterized by the set of weights and biases ${\Theta^{E}}$, and it maps the full-order solution to the kPOD reduced order coefficients:

\begin{equation}
    \boldsymbol{z}^{\text{AE}}_{\text{NN}} =  \psi^{E}(\boldsymbol{u}_{h};{\Theta^{E}}).
\end{equation}
The decoder function $\psi^{D}(\cdot;{\Theta^{D}})$ is the same as in the decoder-based approach, and it maps the reduced order coefficients back to the full-order solution:

\begin{equation}
    \boldsymbol{u}^{\text{AE}}_{\text{NN}} = \psi^{D}(\boldsymbol{z}^{\text{AE}}_{\text{NN}};{\Theta^{D}}).
\end{equation}
During the learning procedure we seek the optimal configuration of the encoder and decoder weights such that we can approximate the forward and backward mapping as follows:

\begin{equation}
    \left\{\begin{aligned}
            &\psi^E(\boldsymbol{u}_h ; \hat{\Theta}^E)\approx F(\boldsymbol{u}_h), \\
            &\psi^D(\boldsymbol{z}^{\text{AE}}_{\text{NN}} ; \hat{\Theta}^D) \approx F^{-1}(\boldsymbol{z}^{\text{AE}}_{\text{NN}}).
        \end{aligned} \right.
\end{equation}
To find the optimal weights $(\hat{\Theta}^{E},\hat{\Theta}^{D})$, we minimize the joint loss function defined as:

\begin{equation}
\label{ae-opt}
(\hat{\Theta}^{E},\hat{\Theta}^{D})=\underset{(\Theta^{E},\Theta^{D})}{\operatorname{argmin}} \ \mathcal{J}^{I_d}(\Theta^{E},\Theta^{D})+ \lambda \mathcal{J}^{\text{kPOD}}(\Theta^{E}).
\end{equation}
The first term, namely $\mathcal{J}^{I_{d}}$, is the reconstruction loss and measures the difference between the input and output of the autoencoder.
On the other hand, $\mathcal{J}^{\text{kPOD}}$ represents the kPOD loss, and measures the discrepancy between the encoded coefficients and the true kPOD coefficients obtained using Equation \eqref{forward-map}.
The hyperparameter $\lambda \in \mathbb{R}_{+}$ plays a crucial role in balancing the contributions of these two terms during training. It allows us to control the trade-off between the reconstruction accuracy and the faithfulness of the reduced order representation.

\subsection{Neural Network Training}
\label{subsec:training}

The neural network regression in our implementation is performed using two types of architectures: feedforward and convolutional layers.
The characteristics of these architectures are presented in Appendix \ref{app:arch}.

To obtain the optimal parameters, we employ a supervised learning paradigm. Specifically, we have $N_{\text{tr}}$ pairs of full-order solutions and their corresponding kPOD coefficients $\{(\boldsymbol{u}_{h}^{i},\boldsymbol{z}^{i})\}_{i=1}^{N_{\text{tr}}}$. During each epoch, we iterate over the training data and calculate the reconstruction loss, which can be chosen either as the Mean Squared Error (MSE) or the Sinkhorn divergence.

When using the MSE measure, the general loss term can be written as:
\begin{equation}
    \mathcal{J}(\Theta) = \frac{1}{N_{\text{tr}}} \sum_{i=1}^{N_{\text{tr}}} \|\boldsymbol{y}^{i} - \boldsymbol{y}_{\text{NN}}^{i}(\Theta)\|^{2},
\end{equation}
where the expressions for the NN output $\boldsymbol{y}_{\text{NN}}(\Theta)$ and the desired output $\boldsymbol{y}$ need to be specialized depending on the type of network under consideration, as shown in Table \ref{tab:loss-function-comparison}.

On the other hand, when using the Sinkhorn divergence, each sample is treated as a point belonging to the support of a probability distribution. Therefore, we have the following expressions for the distribution produced by the network $\nu_{\text{NN}}$ and the desired distribution $\mu$:

\begin{equation}
    \mu = \frac{1}{N_{\text{tr}}} \sum_{i=1}^{N_{\text{tr}}} \delta_{\boldsymbol{y}^{i}}, \quad \nu_{\text{NN}} = \frac{1}{N_{\text{tr}}} \sum_{j=1}^{N_{\text{tr}}} \delta_{\boldsymbol{y}_{\text{NN}}^{j}(\Theta)} \, .
\end{equation}
These distributions allow us to calculate the loss function as the Sinkhorn divergence by means of:

\begin{equation}
    \mathcal{J}(\Theta) = S_{\epsilon}(\mu, \nu_{\text{NN}}).
\end{equation}
In this case, to ensure the correct coupling of the input and output distributions while preserving their order, an initial training phase using the MSE loss for a limited number of epochs is required.
To prevent overfitting, we implement an early stopping criterion, and we set the total number of epochs to $10^3$.

The neural network architectures are implemented using PyTorch \cite{pytorch}, while the Sinkhorn divergence is computed using the geomLoss package \cite{feydy2019interpolating}.
The experiments are conducted on a machine equipped with an NVIDIA GeForce RTX 3080 GPU and an 11th Gen Intel(R) Core(TM) i7-11700 @ 2.50GHz processor running CentOS Stream release 8.

\renewcommand{\arraystretch}{1.5}
\begin{table}[ht]
    \centering
    \begin{tabular}{|c|c|c|}
        \hline
        \cellcolor[HTML]{E5E3E3}Loss Function & \cellcolor[HTML]{E5E3E3}$\boldsymbol{y}$ & \cellcolor[HTML]{E5E3E3}$\boldsymbol{y}_{\text{NN}}$ \\
        \hline
                    $\mathcal{J}^D$ & \ $\boldsymbol{u}_{h}$ \ & \ $\boldsymbol{u}_{\text{NN}}^D$ \ \\
                      $\mathcal{J}^{I_d}$ & $\boldsymbol{u}_{h}$ & $\boldsymbol{u}_{\text{NN}}^{\text{AE}}$ \\
                      $\mathcal{J}^{\text{kPOD}}$ & $\boldsymbol{z}$ & $\boldsymbol{z}_{\text{NN}}^{\text{AE}}$ \\
        \hline
    \end{tabular}
    \caption{Comparison of Loss Functions and Outputs.}
    \label{tab:loss-function-comparison}
\end{table}
\section{Numerical results}%
\label{sec:numerical_results}

In this section, we test the proposed framework on three different benchmarks in model order reduction  for parametrized PDEs: (i) a Poisson equation with a moving source, (ii) the Burgers' equation, and (iii) an advection-dominated problem.
Importantly, all of these problems exhibit slow decay of the Kolmogorov $n$-width, which pose a challenge for efficient model reduction techniques.

In each test case we employ the finite element discretization to obtain the full-order approximation. Furthermore, all the problems are solved on a 2D unitary domain, except for the Burgers' equation where we have one spatial dimension and the time evolution.

We build the original dataset by solving full-order problems for each instance of the parameter in the designated parameter set $\mathcal{P}_{h} = \{\boldsymbol{\mu}_{j}\}_{j = 1}^{N_{s}}$. Next, the set $\mathcal{P}_{h}$ is split into the training set $\mathcal{P}_\text{tr}$ and the testing set $\mathcal{P}_\text{te}$. The cardinalities of these two sets are dictated by the choice of the training rate $r_t$, such that $N_\text{tr}=r_{t}N_{s}$ and $N_\text{te}=N_{s}-N_{\text{tr}}$.

It is important to note that the solutions to the proposed test cases are not probability density functions. Therefore, in order to compute the Wasserstein distance between a pair of snapshots, we normalize the solutions with respect to the total mass. This normalization ensures that the computed Wasserstein distance captures the discrepancy in the spatial distribution of the solutions, rather than the differences in their magnitudes.

As concerns the error analysis, we use as evaluation metric the relative error between the solution $\boldsymbol{u}_{h}(\boldsymbol{\mu})$ and its approximation $\tilde{\boldsymbol{u}}_{\text{h}}$ for a given parameter $\boldsymbol{\mu}$ and its mean over the testing set $\mathcal{P}_{\text{te}}$, namely:
$$
\varepsilon(\boldsymbol{\mu})=\frac{\left\|\boldsymbol{u}_{h}(\boldsymbol{\mu})-\tilde{\boldsymbol{u}}_{h}(\boldsymbol{\mu})\right\|_{L^{2}}}{\left\|\boldsymbol{u}_h(\boldsymbol{\mu})\right\|_{L^{2}}} \quad \text { and } \quad \bar{\varepsilon} = \frac{1}{N_{\text{te}}} \sum_{\boldsymbol{\mu} \in \mathcal{P}_{\text{te}}} \varepsilon (\boldsymbol{\mu}).
$$
Additionally, we will show the relative error field for selected values in the testing set, given by:
$$
\delta \boldsymbol{u}(\boldsymbol{\mu})=\frac{|\boldsymbol{u}_{h}(\boldsymbol{\mu})-\tilde{\boldsymbol{u}}_{h}(\boldsymbol{\mu})|}{\left\|\boldsymbol{u}_h(\boldsymbol{\mu})\right\|_{L^{2}}} \ .
$$
To provide a comprehensive comparison, we compare each test case individually with state-of-the-art linear and nonlinear approaches: POD projection, POD-Galerkin projection, and DL-ROM.

\subsection{Poisson Equation}%
\label{sub:poission_decoder}
We study the solution of a Poisson equation on the 2D unit square $\Omega = [0,1]\times [0,1]$. The equation and boundary conditions can be summarized as follows:

\begin{equation}
    \begin{cases}
        \Delta u = f & \text{in } \Omega, \\
        u = 0 & \text{on } \partial \Omega,
    \end{cases}
\end{equation}
where the source term $f$ is parameterized as a Gaussian function centered at $\boldsymbol{\mu}=(\mu_1, \mu_2)$ with $\sigma = 0.1$:
\begin{equation}
    f = \frac{100}{2\sigma}\exp\left(-\frac{(x - \mu_1)^2 + (y - \mu_2)^2}{2\sigma^2}\right).
\end{equation}
The numerical approximation employs a finite element discretization using $ P1$ elements on a $32 \times 32$ structured mesh, treated as a 2D image with a single color channel, which is essential for the convolutional autoencoder.

For our study, we consider the parameter space $\mathcal{P} = \{(\mu_1, \mu_2) \mid \mu_1, \mu_2 \in [0.2, 0.8]\} \subset \mathbb{R}^{2}$, that we discretize with a uniform grid consisting of $N_{s} = 100$ samples, yielding the its discrete counterpart $\mathcal{P}_h$. This allows us to investigate the influence of different source locations on the scalar field's solution.

To address the forward problem, we collect the dataset by solving the Poisson equation for values in $\mathcal{P}_h$, and then we use the corresponding subset $\mathcal{P}_{\text{tr}}$ with $r_{t}=70\%$ to compute the Gram matrix using the Wasserstein kernel.

\begin{figure}[h]
\includegraphics[width=0.6\textwidth,trim=10 12 10 25, clip]{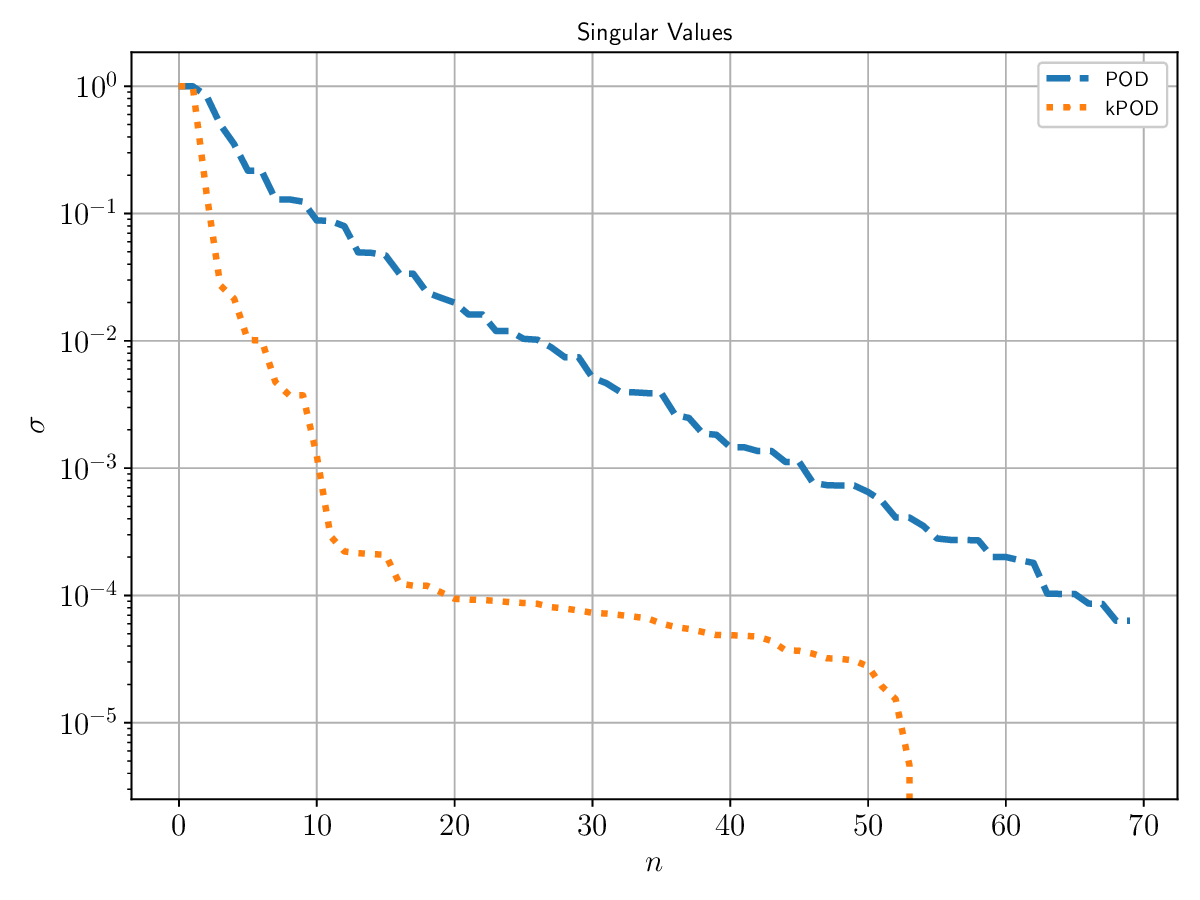}
\caption{Decay of the normalized singular values.}
\label{fig:sing-values-indicators}
\end{figure}

When using kernels based on the Wasserstein distance, the computation of the Gram matrix could represent a computational bottleneck, but given its definition and symmetry we can easily parallelize the computation and consider only its upper triangular part.

Figure \ref{fig:sing-values-indicators} illustrates the comparison between the decay of the singular values for both the POD and the kPOD.

\begin{figure}
\centering
\begin{minipage}{0.48\textwidth}
\includegraphics[width=\textwidth, trim=400 65 410 90, clip]{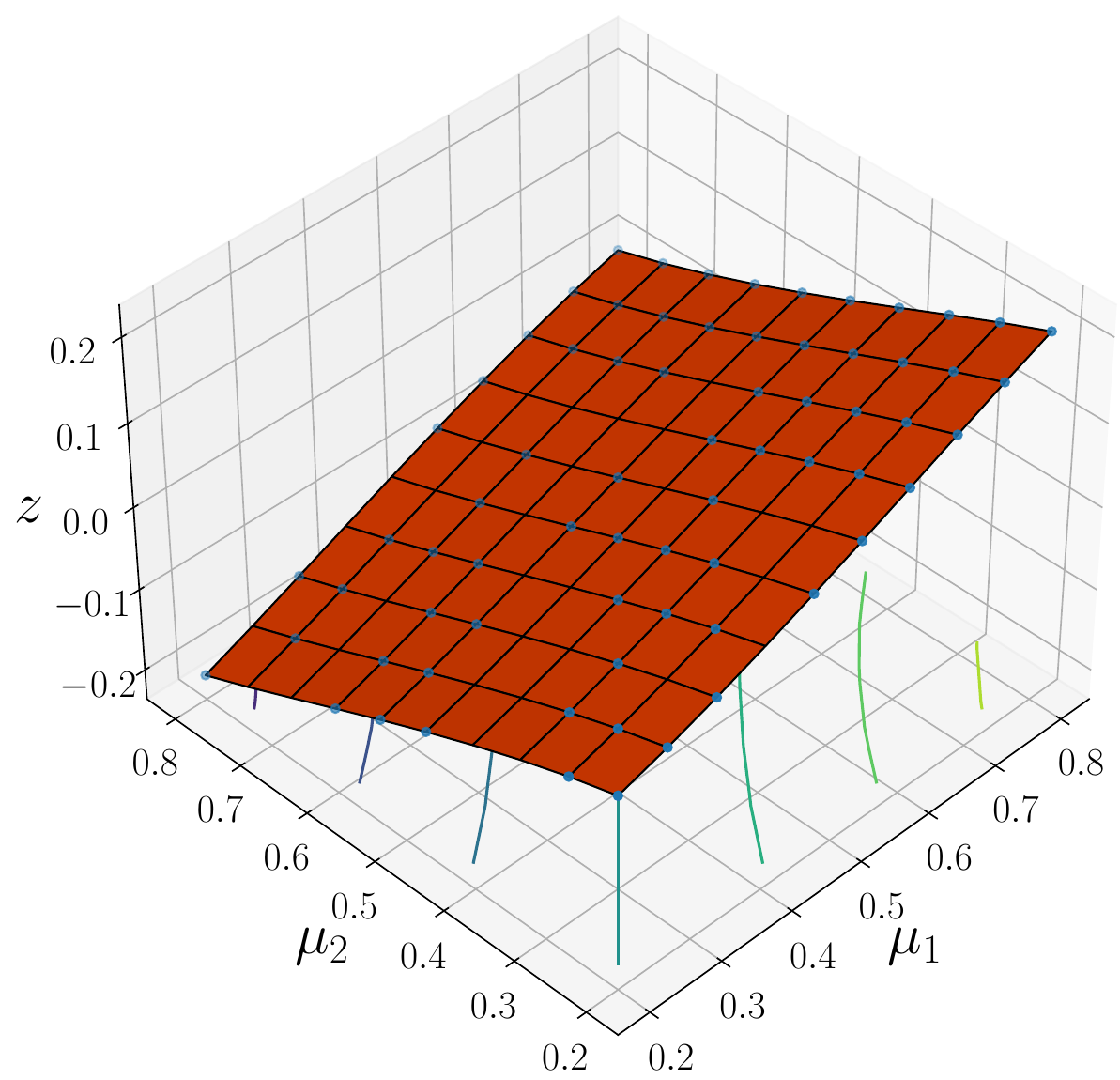}
\end{minipage}
\hfill
\begin{minipage}{0.48\textwidth}
\includegraphics[width=\textwidth, trim=400 65 410 90, clip]{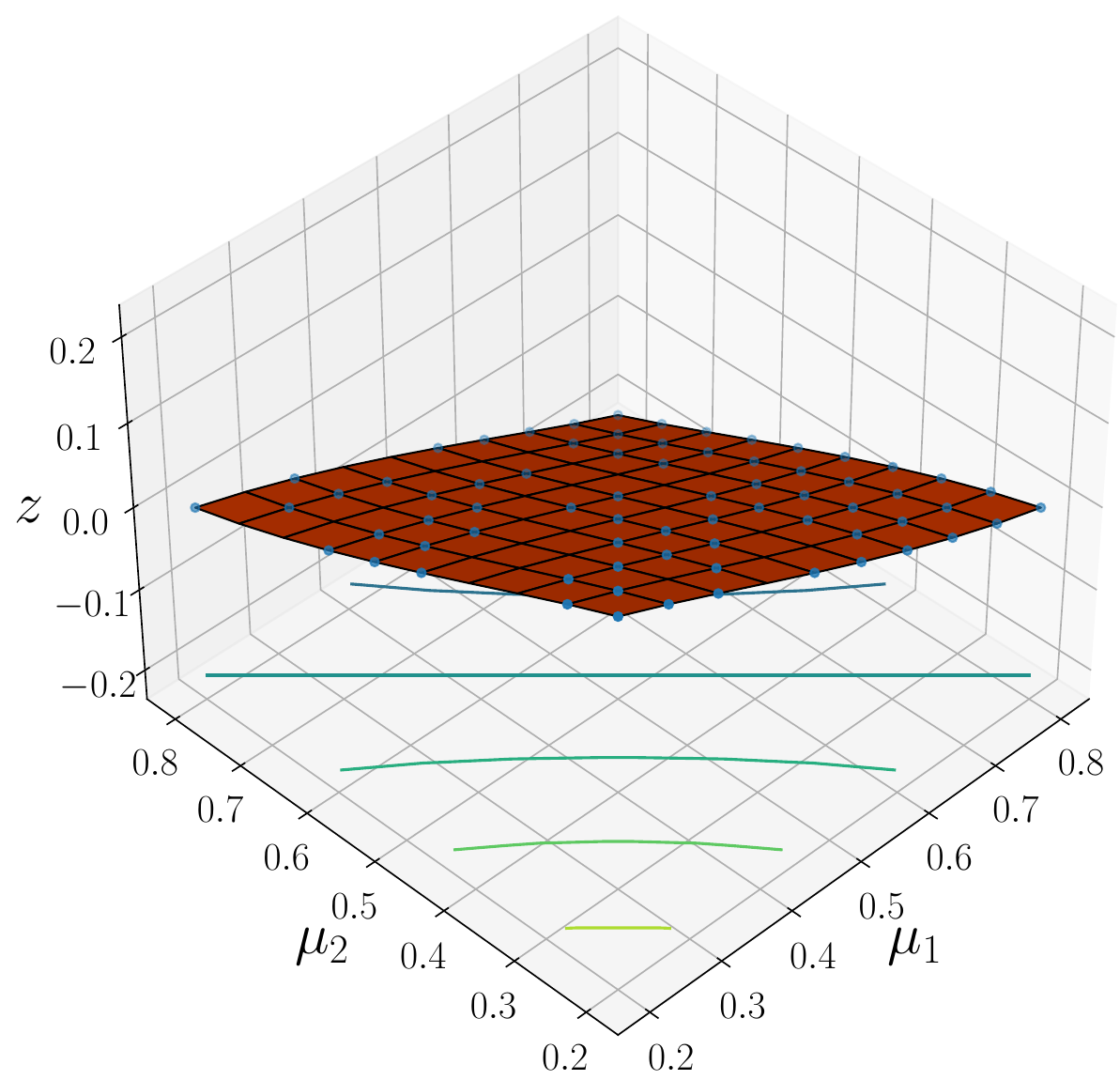}
\end{minipage}
\caption{First (left) and second (right) kPOD coefficients w.r.t.\ the source location. Training points represented in blue.}
\label{fig:kpod_coeffs}
\end{figure}

\begin{figure}
\centering
\begin{minipage}{0.48\textwidth}
\includegraphics[width=\textwidth, trim=400 65 410 90, clip]{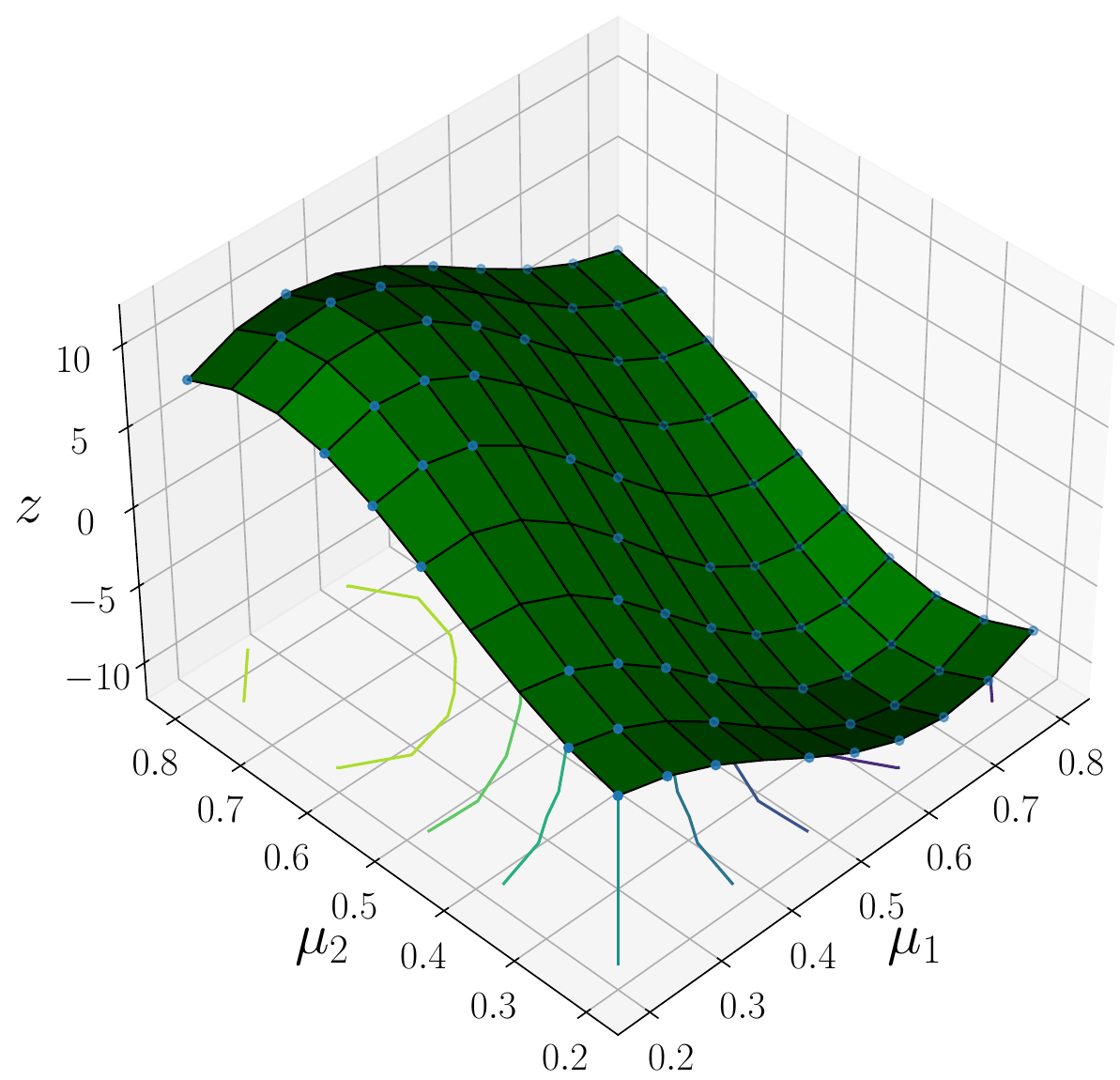}
\end{minipage}
\hfill
\begin{minipage}{0.48\textwidth}
\includegraphics[width=\textwidth, trim=400 65 410 90, clip]{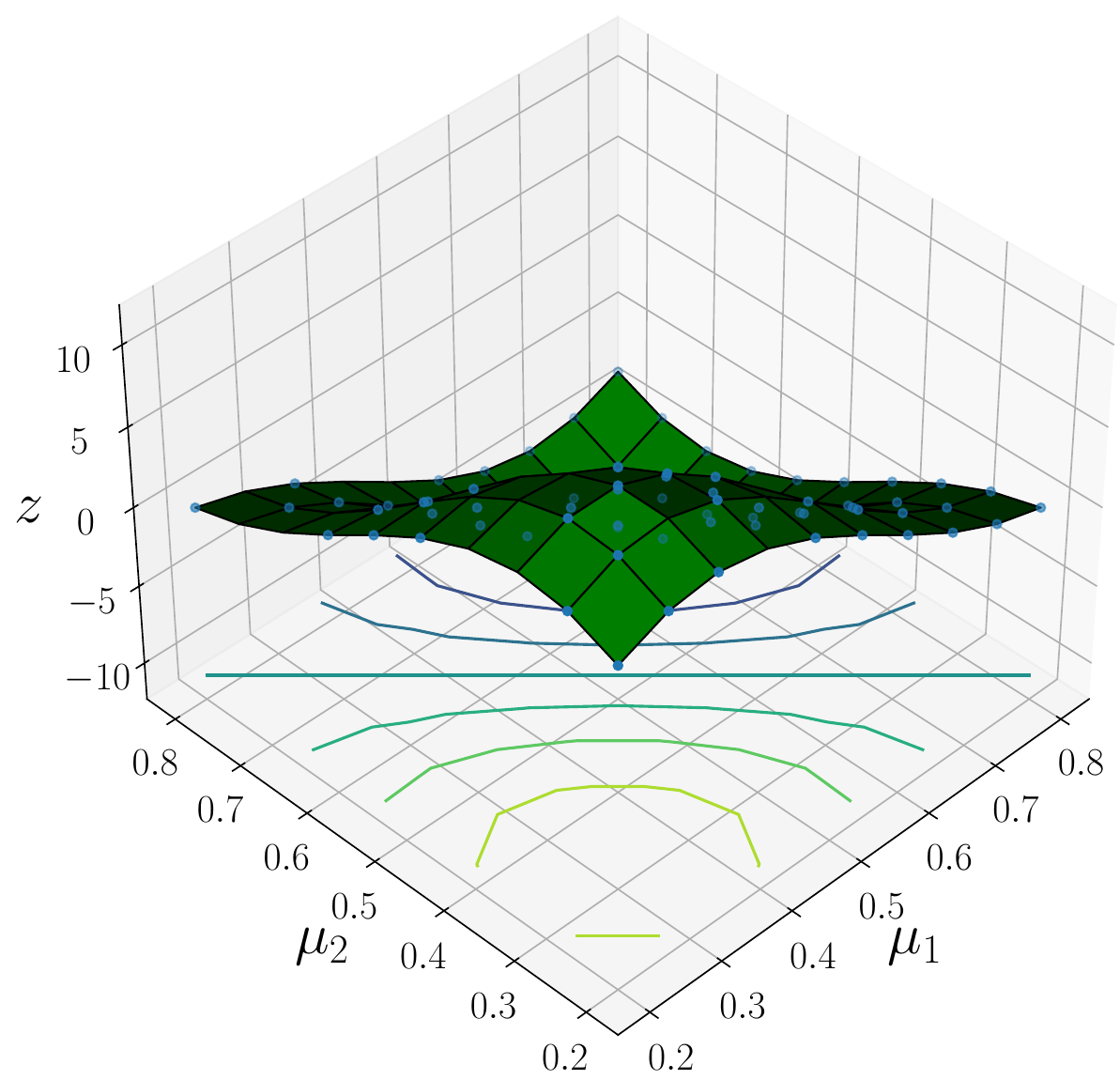}
\end{minipage}
\caption{First (left) and second (right) POD coefficients w.r.t.\ the source location. Training points represented in blue.}
\label{fig:pod_coeffs}
\end{figure}

As expected, the use of the Wasserstein kernel for the kPOD leads to a faster decay of the singular values. We proceed with the projection of the snapshots using the forward map defined in Equation \eqref{forward-map}. In Figure \ref{fig:kpod_coeffs} and \ref{fig:pod_coeffs}, we show the trend of the first two reduced coefficients. We observe that in the case of kPOD, contrarily to what happens within the standard POD approach, the first two reduced coefficients can be accurately approximated by a plane. This implies the existence of a linear relationship between the two parameters under consideration and the kPOD coefficients.

Table \ref{results:poisson} displays the mean relative error $\bar{\varepsilon}$ over the testing set for various models. These models are built upon two architectural approaches: convolutional autoencoder (CAE) and fully connected feedforward autoencoder (FF) (see Appendix \ref{app:arch}), and for both we set the dimension of the reduced representation, \ie\ the number of bottleneck nodes, to $k=5$. The training process employs two different loss functions: mean squared error (MSE) and Sinkhorn loss.

\newcommand{\gradientcell}[6]{
    \ifdimcomp{#1pt}{>}{#3 pt}{#1}{
        \ifdimcomp{#1pt}{<}{#2 pt}{#1}{
            \pgfmathparse{int(round(100*(#1/(#3-#2))-(\minval*(100/(#3-#2)))))}
            \xdef\tempa{\pgfmathresult}
            \cellcolor{#5!\tempa!#4!#6} #1
    }}
}

\newcommand*{\opacity}{40}
\newcommand*{\minval}{0.50}
\newcommand*{\maxval}{2.84}
\begin{table}[htbp]
    \centering
    \begin{tabular}{|c|c|c|c|c|c|}
        \hhline{~|~|----}
        \multicolumn{2}{c|}{} & \multicolumn{2}{c|}{\cellcolor[HTML]{E5E3E3}\textbf{kPOD}} & \multicolumn{2}{c|}{\cellcolor[HTML]{E5E3E3}\textbf{POD}} \\
        \cline{2-6}
        \multicolumn{1}{c|}{} & \diagbox[height=3em, width=12 em]{Loss type}{Architecture} &{FF} & {CAE} & {FF} & {CAE} \\
        \hline
        \multirow{2}{*}{\cellcolor[HTML]{E5E3E3}} & Sinkhorn & \gradientcell{2.02}{\minval}{\maxval}{low}{high}{\opacity} \%& \gradientcell{0.50}{\minval}{\maxval}{low}{high}{\opacity} \%& \gradientcell{1.95}{\minval}{\maxval}{low}{high}{\opacity} \%& \gradientcell{1.00}{\minval}{\maxval}{low}{high}{\opacity} \%\\
        \cellcolor[HTML]{E5E3E3}\raisebox{1.6ex}[1.6ex]{\textbf{Autoencoder}} & MSE & \gradientcell{2.22}{\minval}{\maxval}{low}{high}{\opacity} \%& \gradientcell{0.70}{\minval}{\maxval}{low}{high}{\opacity} \%& \gradientcell{2.08}{\minval}{\maxval}{low}{high}{\opacity} \%& \gradientcell{1.20}{\minval}{\maxval}{low}{high}{\opacity} \%\\
                                                  \hline
        \multirow{2}{*}{\cellcolor[HTML]{E5E3E3}} & Sinkhorn & \gradientcell{2.60}{\minval}{\maxval}{low}{high}{\opacity} \% & \gradientcell{1.50}{\minval}{\maxval}{low}{high}{\opacity} \%& \gradientcell{2.84}{\minval}{\maxval}{low}{high}{\opacity} \%& \gradientcell{1.75}{\minval}{\maxval}{low}{high}{\opacity} \%\\
        \cellcolor[HTML]{E5E3E3}\raisebox{1.6ex}[1.6ex]{\textbf{Decoder}}& MSE & \gradientcell{2.22}{\minval}{2.60}{low}{high}{\opacity} \%& \gradientcell{2.42}{\minval}{\maxval}{low}{high}{\opacity} \%& \gradientcell{2.80}{\minval}{\maxval}{low}{high}{\opacity} \%& \gradientcell{2.61}{\minval}{\maxval}{low}{high}{\opacity} \%\\
        \hline
    \end{tabular}
    \caption{Comparison of mean relative errors for different neural network architectures and loss functions. The models tested include a convolutional autoencoder (CAE) and a fully connected feedforward autoencoder (FF). POD and kPOD refer to the reduced approach exploited to obtain the latent representation.}
    \label{results:poisson}
\end{table}

In our evaluation, we first examined the performance of both autoencoder and decoder approaches. Interestingly, the decoder-only approach consistently exhibited higher reconstruction errors compared to the autoencoder across all models, suggesting that leveraging the full autoencoder architecture is crucial for achieving better reconstruction results. Consequently, we will focus our discussion on the autoencoder-based models, which demonstrate more promising performance w.r.t.\ decoder-only ones (Table \ref{results:poisson}).

\begin{figure}[h]
    \centering
    \includegraphics[width=\linewidth]{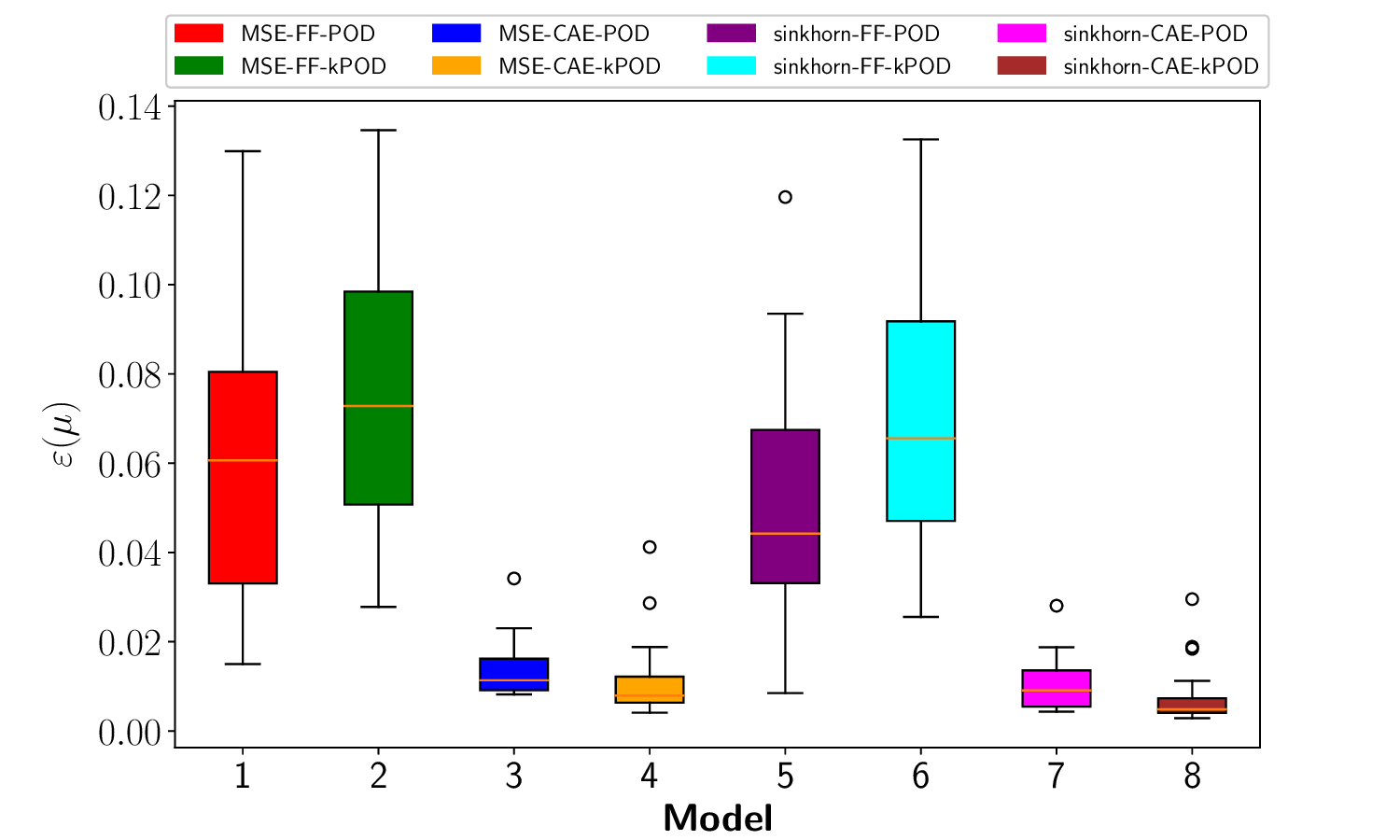}
    \caption{Box-plot of the relative errors for the Poisson test case and the autoencoder architecture.}%
    \label{fig:box_plot_poisson}
\end{figure}

As it can be seen from Figure \ref{fig:box_plot_poisson}, comparing the models using the MSE loss, we observe that the CAE models exhibit lower reconstruction errors compared to their FF counterparts. Specifically, the CAE-kPOD model achieves the lowest reconstruction error among all the MSE-based models, both for the autoencoder and decoder approaches. This result suggests that incorporating the kPOD technique into the CAE architecture enhances its ability to capture the underlying features and achieve more accurate reconstructions.

When considering the sinkhorn loss, we observed that it consistently led to lower reconstruction errors for both the FF and CAE architectures. The CAE-kPOD model, in particular, achieves the lowest reconstruction error among all models, irrespective of the autoencoder or decoder configuration. This result highlights the effectiveness of the sinkhorn loss function in conjunction with the CAE-kPOD model, suggesting that it provides valuable regularization benefits that contribute to improved reconstruction accuracy.

\begin{figure}[h]
\centering
\begin{subfigure}[t]{0.45\linewidth}
\centering
\includegraphics[width=\linewidth, trim=100 0 290 80, clip]{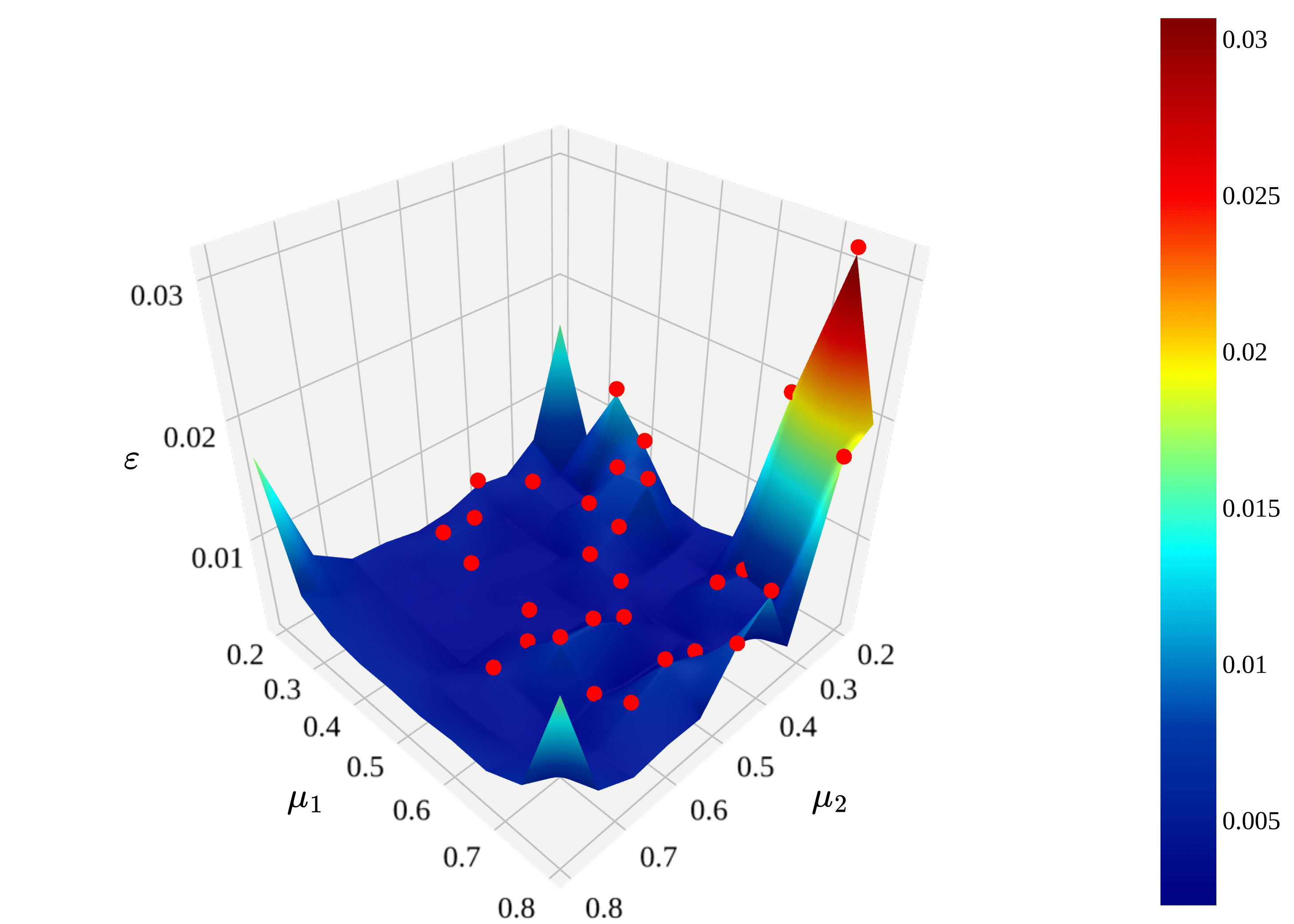}
\caption{Sinkhorn-CAE-kPOD: the tesing points are represented in red.}
\label{fig:err-sink-cae-kpod-poisson}
\end{subfigure}\qquad
\begin{subfigure}[t]{0.45\linewidth}
\centering
\includegraphics[width=\linewidth, trim=100 0 290 80, clip]{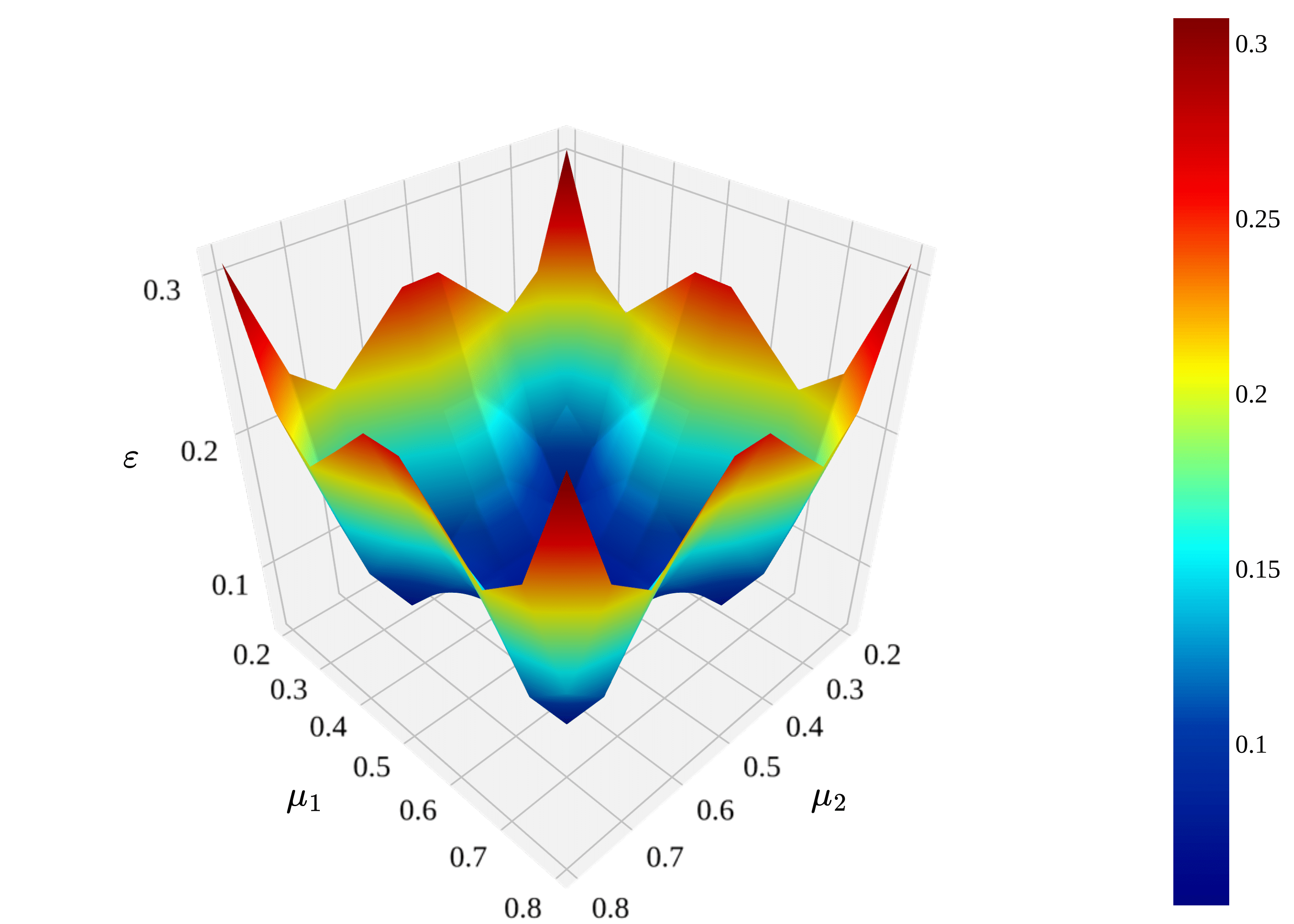}
\caption{Projection on the POD basis with $k = 5$.}
\label{fig:projection-error}
\end{subfigure}
\caption{Relative error with respect to the source location.}
\label{fig:error-comparison}
\end{figure}

\begin{figure}[h]
\centering
\begin{subfigure}[t]{1\linewidth}
    \caption{Solution fields}
    \includegraphics[width=1\linewidth,trim=100 0 50 0, clip]{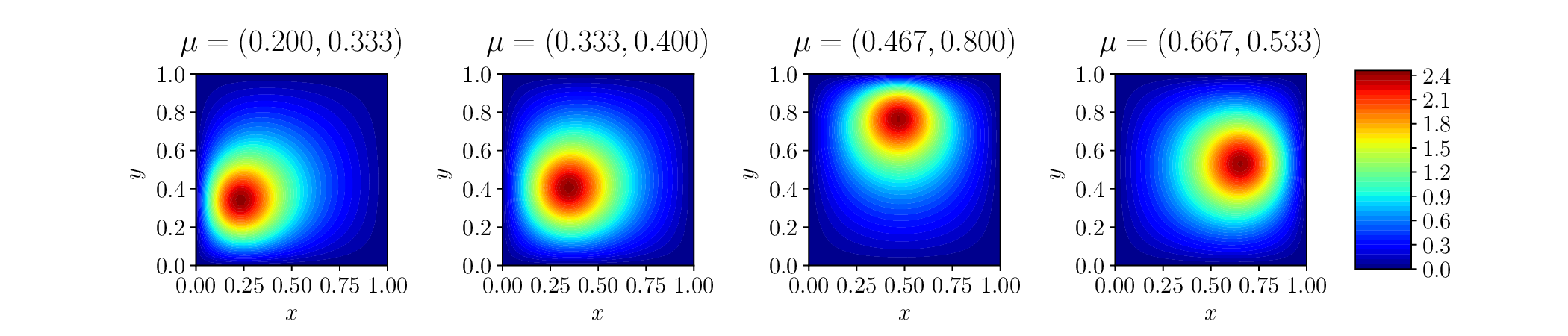}
    \label{fig:fom_poisson}
\end{subfigure}
\begin{subfigure}[t]{1\linewidth}
    \centering
    \caption{Relative error fields}
    \includegraphics[width=1\linewidth,trim=100 0 50 0, clip]{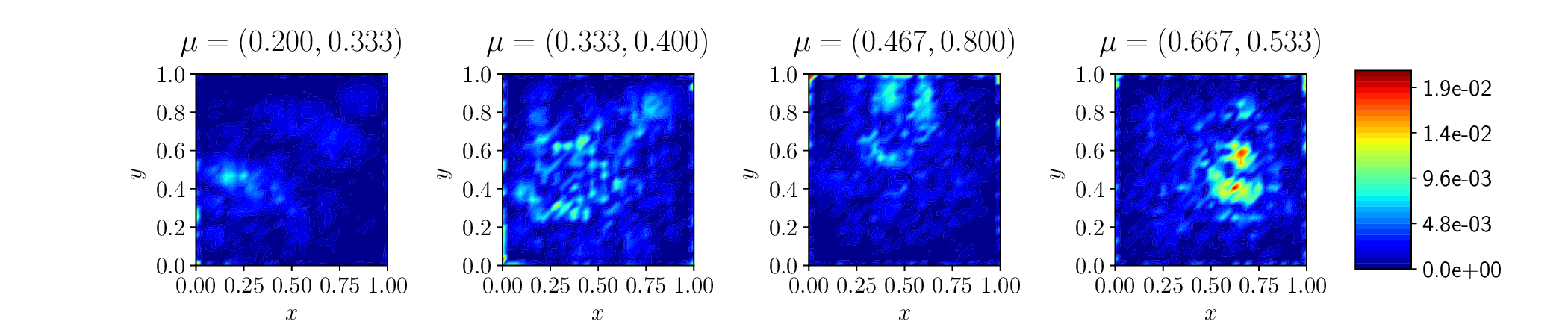}
    \label{fig:error_poisson}
\end{subfigure}
\caption{Solution and relative error fields for the Poisson equation for different positions of the source term.}
\label{fig:poisson}
\end{figure}

The variation of the relative error as a function of the source position for the sinkhorn-CAE-kPOD model is shown in Figure \ref{fig:err-sink-cae-kpod-poisson}, whereas Figure \ref{fig:poisson} presents representative solutions and their corresponding relative error fields for parameters in $\mathcal{P}_{\text{te}}$. It can be observed that the error does not exceed 3\% and takes on higher values when the testing points correspond to parameter values near the boundary of the parameter space, as expected in such cases.

The effectiveness of our approach becomes evident when compared to the error of direct projection onto the POD basis (see Table \ref{table:projection-error}). By directly projecting the FOM solutions onto a basis with the same cardinality as the bottleneck of our autoencoder, \ie \ $k=5$, the error is significantly higher, with an average of $14.66\%$ (see Figure \ref{fig:projection-error}). As the number of basis increases, the error gradually decreases, but it only becomes lower than the sinkhorn-CAE-kPOD model when $k = 30$.

However, it is essential to acknowledge that the direct projection error represents only one component of the overall error for linear methods. Therefore, when considering a generic linear method with the same number of bases, we can reasonably expect even higher errors compared to the performance achieved by the kPOD based approach.

\begin{table}[htbp]
\centering
\begin{tabular}{cc}
\toprule
Number of basis  & Projection relative error (\%) \\
\midrule
5  & 14.66 \\
10 & 5.20 \\
20 & 1.19 \\
30 & 0.34 \\
\bottomrule
\end{tabular}
\caption{Mean relative errors for different numbers of basis functions.}
\label{table:projection-error}
\end{table}

\subsection{Advection Dominated Problem}%
\label{sub:advection_dominated}
Here, we are interested in approximating the solution of a parametrized advection-diffusion problem defined on the unit square $\Omega = [0,1]\times [0,1]$:

\begin{equation}
    \label{eq:adv}
    \begin{cases}
        - \alpha(\mu)\Delta u  + \beta \cdot \nabla u = 0 & \text{in } \Omega, \\
        u = 0 & \text{on } \partial \Omega,
    \end{cases}
\end{equation}
where $\beta$ and $\alpha$ represent the advection and diffusion terms, respectively, and homogeneous Dirichlet boundary conditions are imposed on the boundary.
We assume that the convective field is constant throughout the domain, and given by $\beta = (1,1)$. The diffusivity is parameterized by $\mu \in [0, 10]$ by defining $\alpha(\mu) = 10^{-\mu}$, which is inversely proportional to the Péclet number of the model.

The offline procedure begins by creating the discrete set $\mathcal{P}_{h}$, comprising samples drawn uniformly from the parameter set $\mathcal{P}$. We use $N_s = 100$ samples and set the training ratio to $r_t = 70\%$.

\begin{figure}[h]
\centering
\begin{subfigure}[t]{1\linewidth}
    \caption{Solution fields}
    \includegraphics[width=1\linewidth,trim=70 0 50 0, clip]{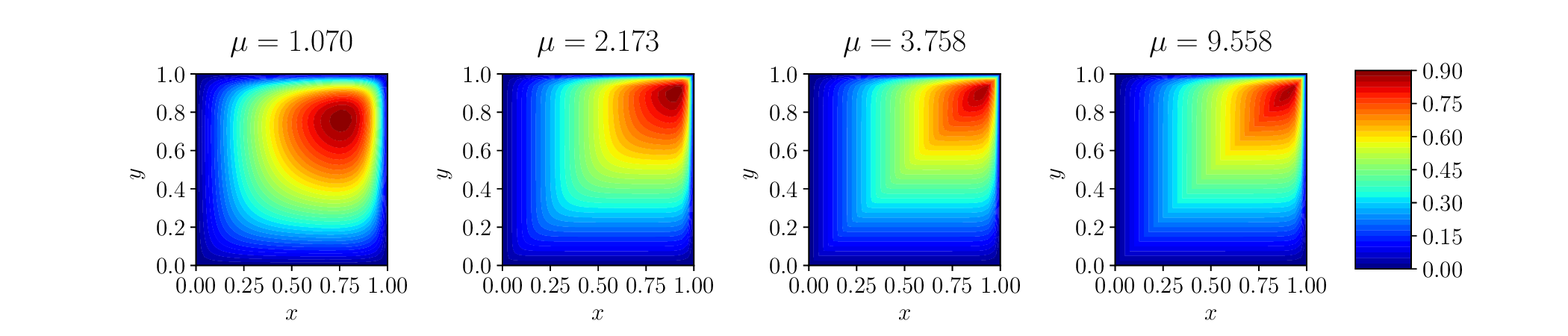}
    \label{fig:fom_advection}
\end{subfigure}
\begin{subfigure}[t]{1\linewidth}
    \centering
    \caption{Relative error fields}
    \includegraphics[width=1\linewidth,trim=70 0 50 0, clip]{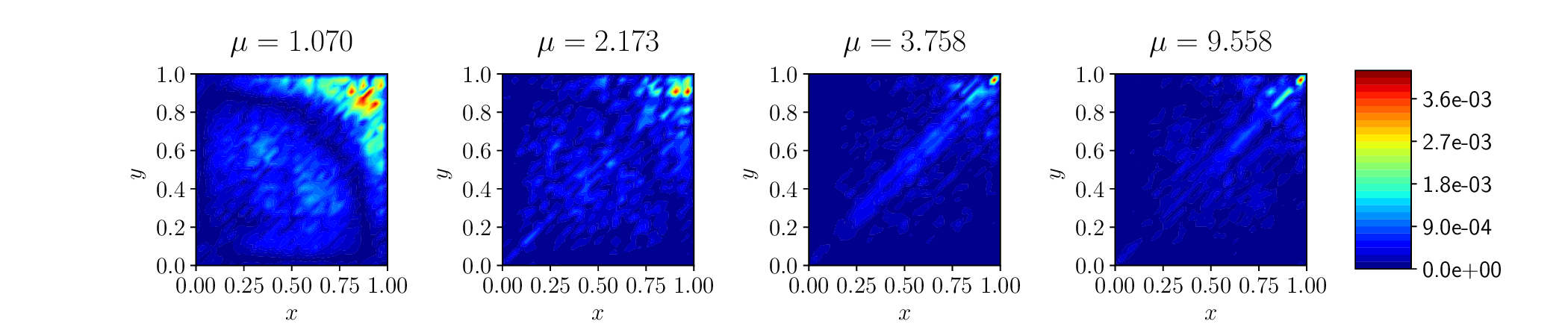}
    \label{fig:error_advection}
\end{subfigure}
\caption{Solutions and relative errors for the advection-dominated problem with varying diffusivity values.}
\label{fig:advection}
\end{figure}

Figure \ref{fig:advection} presents some representative solutions and their corresponding relative error fields, obtained with the Sinkhorn-CAE-kPOD model. It can be observed that, when $\mu$ increases, the convective contribution dominates the diffusive term, resulting in significant displacement of the solution in the direction of $\beta$.

Table \ref{table:results-advection} presents the mean reconstruction relative errors using the autoencoder approach, w.r.t.\ the test set $\mathcal{P}_{\text{te}}$.

Firstly, it is observed that for this particular problem, when using the MSE loss function, the fully connected architecture outperforms the convolutional autoencoder. This is probably due to the thin boundary layer developed in the top-right corner, which of cours affects local convolutions. Moreover, this highlights the importance, and yet the issue, of choosing an architecture that is specifically suited to the problem at hand.
Despite this, utilizing the Sinkhorn Loss during training consistently leads to improved model accuracy, and the lowest error is achieved by the CAE architecture, once again highlighting the effectiveness of embdedding OT-based loss in the proposed neural network framework.

\renewcommand*{\minval}{0.069}
\renewcommand*{\maxval}{0.189}
\begin{table}[htbp]
\centering
\begin{tabular}{|c|c|c|}
\hhline{~|--}
\multicolumn{1}{c|}{} & \cellcolor[HTML]{E5E3E3}\textbf{CAE} & \cellcolor[HTML]{E5E3E3}\textbf{FF} \\
\hline
\cellcolor[HTML]{E5E3E3}\textbf{Sinkhorn} & \gradientcell{0.069}{\minval}{0.189}{low}{high}{\opacity} \% & \gradientcell{0.081}{\minval}{0.189}{low}{high}{\opacity} \% \\
\hline
\cellcolor[HTML]{E5E3E3}\textbf{MSE}& \gradientcell{0.189}{\minval}{0.189}{low}{high}{\opacity} \% & \gradientcell{0.127}{\minval}{0.189}{low}{high}{\opacity} \% \\
\hline
\end{tabular}
\caption{Summary of results for different loss functions and architectures.}
\label{table:results-advection}
\end{table}

In the following, we aim at comparing these results with those obtained through the POD-Galerkin (POD-G) reduced order method, which involves projecting the system of equations \eqref{eq:adv} onto the POD modes\footnote{To obtain physically meaningful solutions, one needs to stabilize the system using the SUPG method \cite{supg-pod-rom}.}.

\begin{figure}
\centering
\begin{minipage}{0.49\textwidth}
\includegraphics[width=\textwidth,trim=10 10 10 10, clip]{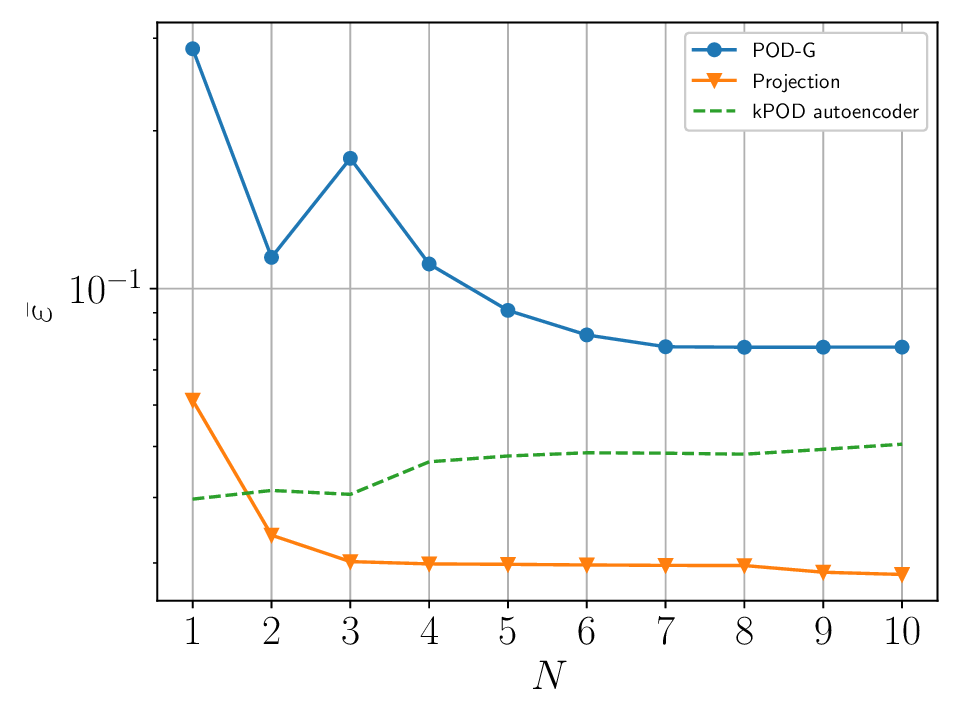}
\put(-140,160){$N_{\text{tr}}=10;N_{\text{te}}=90$}
\end{minipage}
\hfill
\begin{minipage}{0.49\textwidth}
\includegraphics[width=\textwidth,trim=10 10 10 10, clip]{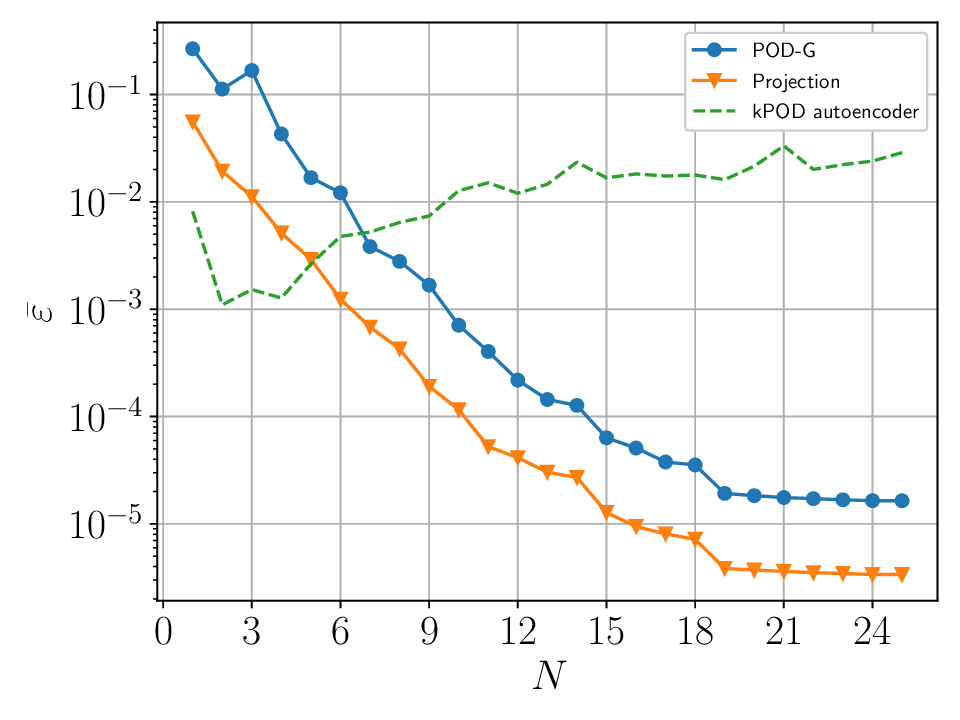}
\put(-140,160){$N_{\text{tr}}=50;N_{\text{te}}=50$}
\end{minipage}
\caption{Mean relative error with $N_{\text{tr}}$=10 (left) and $N_{\text{tr}}$=50 (right).}
\label{fig:error_POD_G_adv}
\end{figure}

Figure \ref{fig:error_POD_G_adv} depicts the mean relative errors obtained on the testing set for the POD-G method and the kPOD autoencoder, while varying the number of basis functions or bottleneck nodes used.
The left plot showcases the errors when using $r_t = 10\%$ of the dataset for training, whereas the right plot corresponds to $r_t = 50\%$.

Despite the linearity of the problem greatly enhances the effectiveness of the POD-G method, a clear limitation arises when utilizing a small number of snapshots for constructing the reduced basis. In such a context, the accuracy of the model saturates for a limited number of bases and deteriorates upon the addition of new bases. In this scenario, the autoencoder exhibits better reconstruction error, although still inferior to direct projection onto POD modes.
It is noteworthy that when 50\% of the samples are used for training, the situation is markedly different. The error for POD-G decreases rapidly as the cardinality of the reduced basis expands.  It is crucial to note that error saturation now occurs for  $N > 19$, and that the autoencoder can obtain better accuracy for $N<6$, before entering the overfitting regime, showing greater compression capabilities.

\subsection{Viscous Burgers Equation}%
\label{sub:burgers_equation}
As a final test case, we consider the 1-dimensional Burgers equation:

\begin{equation}
\begin{cases}
\frac{\partial u}{\partial t} + u \frac{\partial u}{\partial x} = \nu \frac{\partial^2 u}{\partial x^2}, \quad \forall x \in [0,1],\\
u(0,t) = u(1,t) = 0, \quad \forall t > 0, \\
u(x,0) = u_0(x, \mu), \quad  \forall x \in [0,1],
\end{cases}
\end{equation}
where $\nu$ represents the viscosity parameter, and the initial condition is chosen as the parametric Heaviside function given by:

\begin{equation}
u_0(x, \mu) = \begin{cases}
0, & \text{if } x < \mu, \\
1, & \text{if } \mu \leq x < 0.2 + \mu, \\
0, & \text{if } x \geq 0.2 + \mu,
\end{cases}
\end{equation}
where $\mu \in  \mathcal{P} = [0.2,0.4]$ represents the shift parameter that determines the location of the Heaviside function within the domain.

The full-order discretization of the problem consists of both spatial and temporal discretizations. For the spatial discretization, we employ the FE method with $P1$ elements. As for the temporal discretization, we utilize the Runge-Kutta 2/3 method (RK23) \cite{rk23}, which provides an accurate time integration scheme.
Although the full-order solution is obtained on a $129\times 129$ spatial/temporal grid, we perform subsampling by selecting equispaced nodes to obtain a reduced image size of $32 \times 32$. This approach allows us to retain the essential features of the solution while reducing the computational cost of training the model.

\begin{figure}[h]
\centering
\begin{subfigure}[t]{1\linewidth}
    \caption{Solution fields}
    \includegraphics[width=1\linewidth,trim=70 0 50 0, clip]{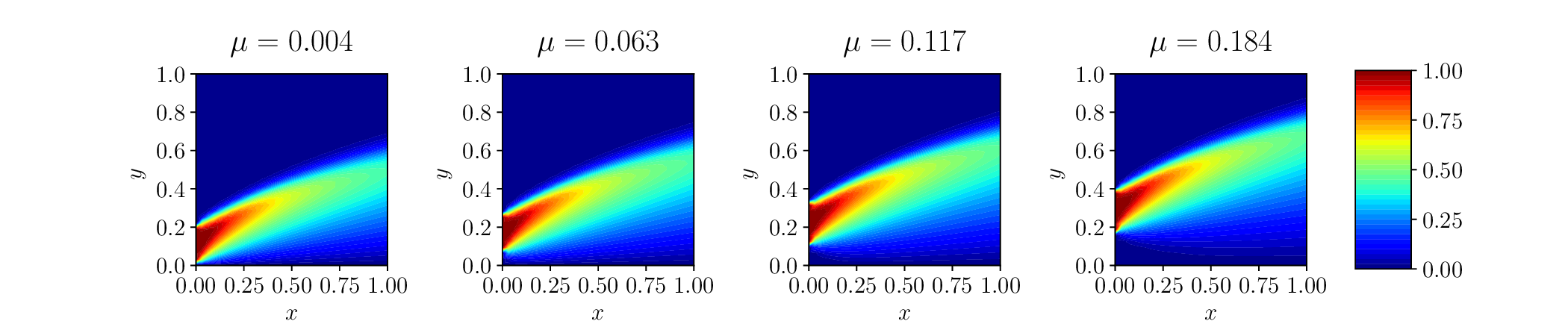}
    \label{fig:fom_advection}
\end{subfigure}
\begin{subfigure}[t]{1\linewidth}
    \centering
    \caption{Relative error fields}
    \includegraphics[width=1\linewidth,trim=70 0 50 0, clip]{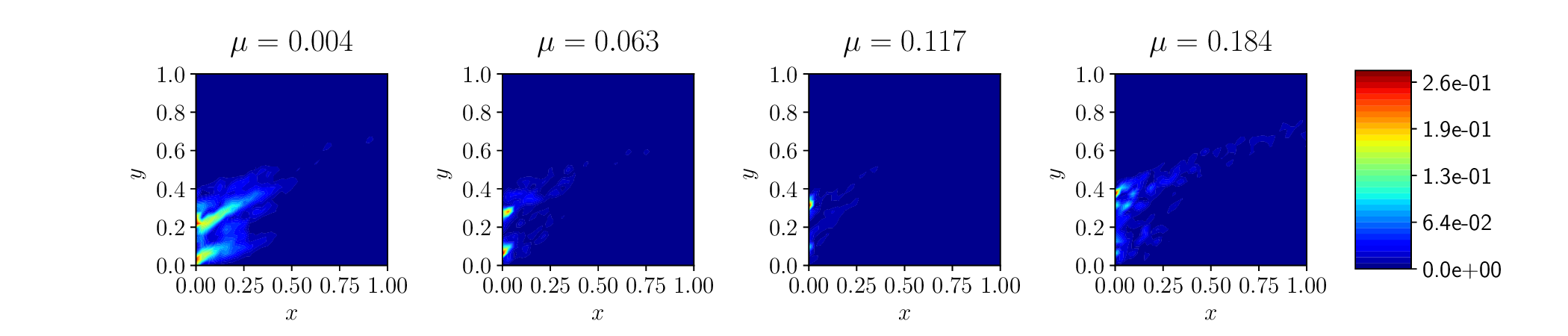}
    \label{fig:error_advection}
\end{subfigure}
\caption{Solutions and relative errors of the Burgers' equation for different initial condition locations.}
\label{fig:burgers}
\end{figure}

As in the previous cases, we proceed with the discretization of the parametric interval considering a grid of $N_s = 100$ equispaced values $\mathcal{P}_{h}$, and select as training rate $r_t = 50\%$.
Figure \ref{fig:burgers} depicts representative solutions and their corresponding relative error fields, obtained from the CAE-kPOD model, while varying the shift parameter.

In this test case, we compare our proposed kPOD approach with the original DL-ROM architecture \cite{fresca2021}. Since the focus is not the optimization of the hyperparameter configuration, we restrict ourselves to using a convolutional autoencoder with the same architecture as our kPOD model. Specifically, we ensure that the autoencoder has the same bottleneck dimension ($k=5$) as our kPOD model but is trained without imposing any specific reduced representation. Moreover the kPOD-CAE is trained using the Sinkhorn loss, whereas the DL-ROM approach is trained using standard MSE.

\begin{figure}[h]
    \centering
    \includegraphics[width=1\textwidth,trim=120 0 120 0, clip]{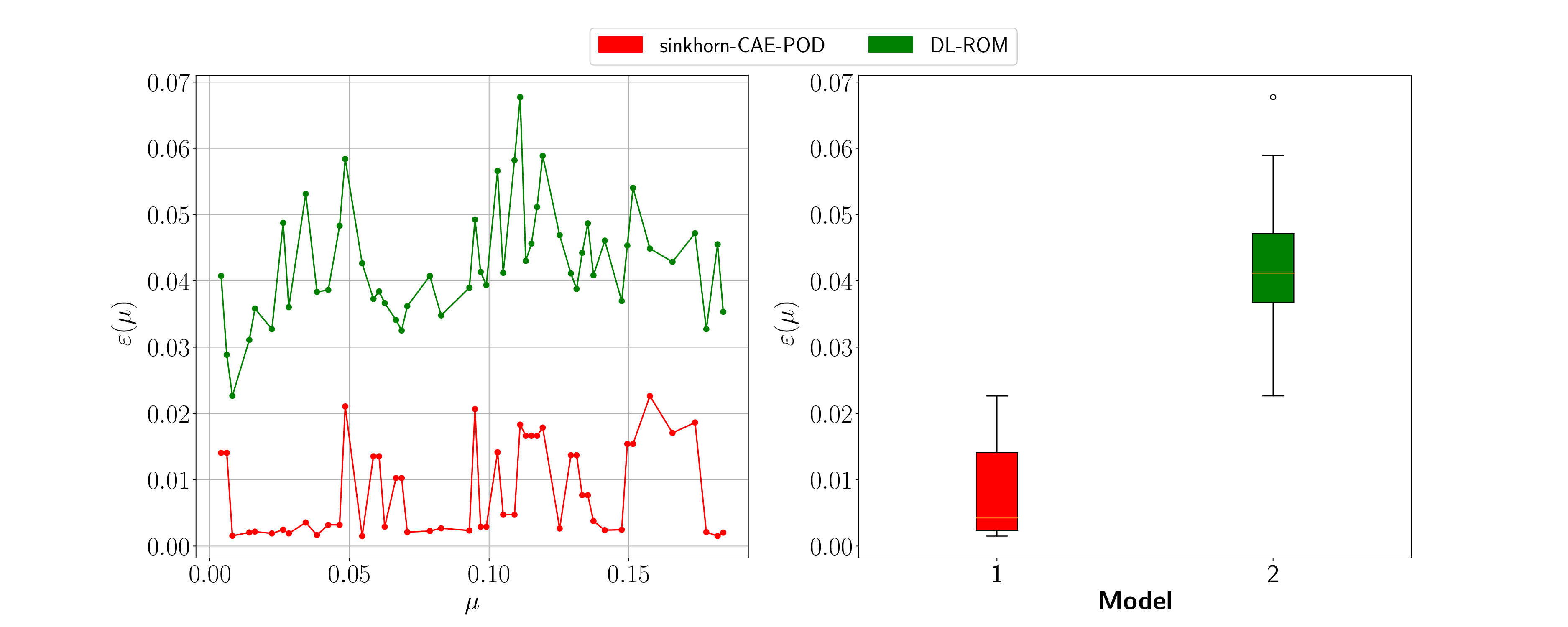}
    \caption{Relative errors for the DL-ROM architecture and the kPOD-autoencoder across the testing set.}%
    \label{fig:error_burgers}
\end{figure}

The reconstruction relative error across $\mathcal{P}_{\text{te}}$ is shown in Figure \ref{fig:error_burgers}. It can be observed that the kPOD-CAE trained with Sinkhorn regularization consistently outperforms its DL-ROM counterpart. This suggests that the reduced representation provided by kPOD is more efficient compared to the one learned through unconstrained training of the autoencoder. This efficiency improvement is reflected in a lower mean error, with $\overline{\varepsilon}_{\text{kPOD}}=8 \times 10^{-3}$ for the kPOD-CAE and $\overline{\varepsilon}_{\text{DL-ROM}}=4 \times 10^{-2}$ for the DL-ROM method.

\subsection{Inviscid Burgers Equation}%
\label{sub:burgers_equation_inviscid}

As an additional challenge for the proposed methodology, we analyze the 1-dimensional Burgers equation without viscosity, commonly exploited to describe shock waves and traffic phenomena:

\begin{equation}
\label{eq:burgers-inviscid}
\begin{cases}
\frac{\partial u}{\partial t} + u \frac{\partial u}{\partial x} = 0, &  x \in [-L, L], t \in [0, T],\\
u(-L, t) = 0, &  t \in [0, T], \\
u(x,0) = u_0(x), &  x \in (-L, L],
\end{cases}
\end{equation}
\begin{equation*}
\text{where} \quad u_0(x) = \begin{cases}
u_-, & \text{if } x \in (-L,0), \\
u_+, & \text{if } x \in [0,L].
\end{cases}
\end{equation*}
The discontinuity values \(u_-\) and \(u_+\) generate a shock wave with speed $s = (u_{+}+u_{-})/2$ given by the Rankine-Hugoniot condition moving through the domain.
Moreover, the inflow condition \(u(-L, t) = 0\) leads to a rarefaction wave originating at the point \((x,t) = (-L,0)\), and affecting the dynamics of the flow.

To capture the evolution of the solution numerically, we use a finite difference scheme. The computational domain is discretized into 200 spatial points over the interval \([-1, 1]\) and evolved over 200 time steps with a step size \(\Delta t = 0.01\). The numerical method involves a first-order upwind scheme for the advection term combined with a high-resolution flux limiter to effectively manage numerical dissipation and preserve the discontinuity at the shock. Time integration is performed using the Runge-Kutta 4/5 method (RK45).

\begin{figure}[h]
\centering
\begin{subfigure}[t]{1\linewidth}
    \includegraphics[width=1\linewidth,trim=70 0 50 0, clip]{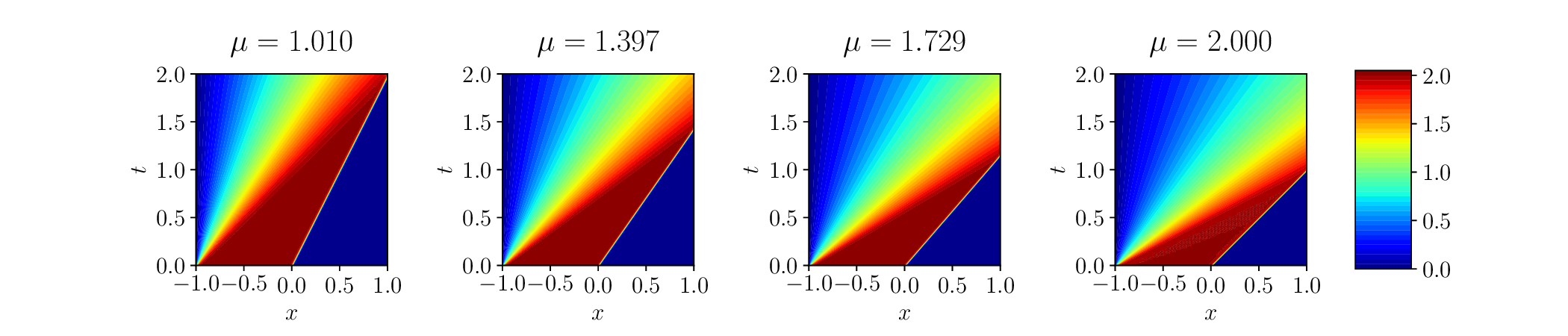}
    \label{fig:fom_burgers_inviscid}
\end{subfigure}\vspace{-1em}
\begin{subfigure}[t]{1\linewidth}
    \centering
    \includegraphics[width=1\linewidth,trim=70 0 50 0, clip]{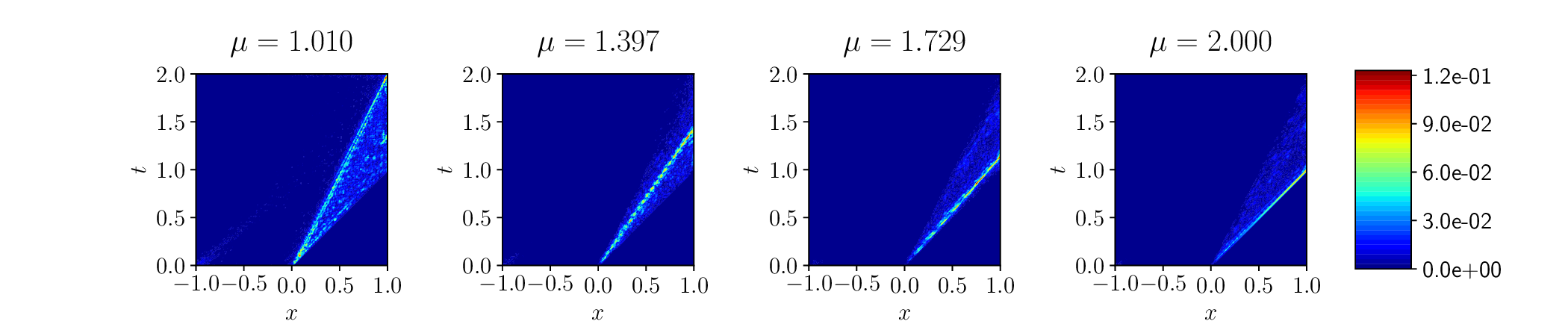}
    \label{fig:error_burgers_inviscid}
\end{subfigure}\vspace{-1.5em}
\caption{Solutions (top row) and relative errors (bottom row) for the inviscid Burgers' equation across various discontinuity values in initial conditions.}
\label{fig:burgers_inviscid}
\end{figure}

We have chosen to parametrize the initial condition by considering \(\mu = u_-
\in \mathcal{P} = [1,2]\) and \( u_+ = 0 \). This parametric variation directly influences the rarefaction wave and the slope of the shock wave, consistently with the Rankine-Hugoniot condition (see Figure \ref{fig:burgers_inviscid}).

We considered \(N_s = 200\) snapshots obtained by solving the problem, with the
parameter \(\mu\) varied on an equispaced grid $\mathcal{P}_{h}$  in the interval \([1,2]\). We used 70\% of the dataset to train a Sinkhorn-CAE-kPOD model with a bottleneck dimension \(k = 5\). Within this test case, we also show the online approximation results for the ROM strategy. Thus, we focused on the task of reconstructing the map from parameters to kPOD coefficients \(\mathscr{F}_{\text{online}}: \mu \in \mathcal{P} \rightarrow \boldsymbol{z}(\mu) \in \mathbb{R}^{k}\).

\begin{figure}[h]
    \centering
    \includegraphics[width=1\linewidth]{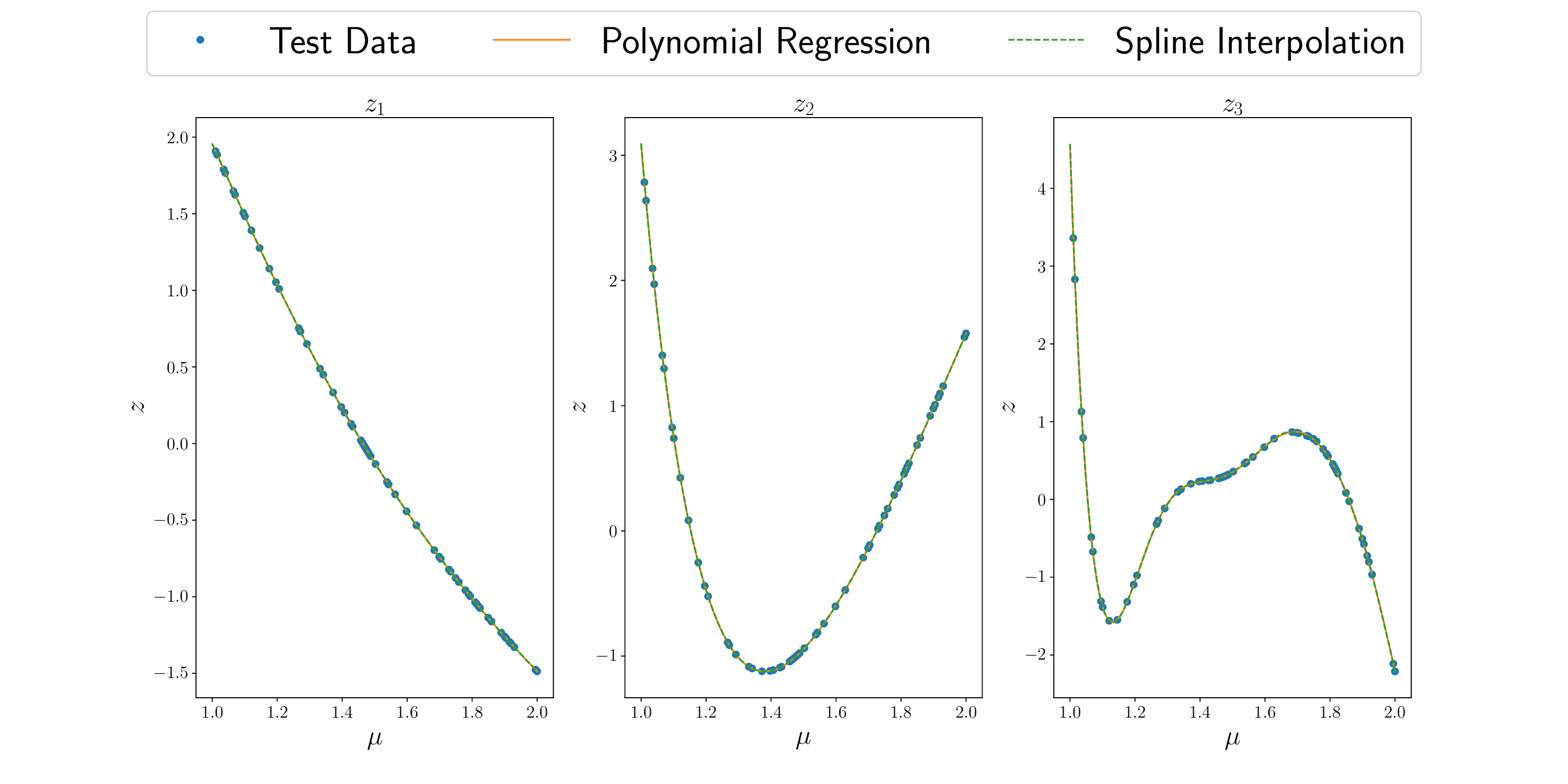}
    \caption{Comparison of polynomial regression with degree 25 and cubic spline interpolation methods showing their performance in predicting the reduced coefficients \( \boldsymbol{z}(\mu) \).}
    \label{fig:polynomial_vs_spline}
\end{figure}

As detailed previously, the use of the kernel within kPOD results in a nonlinear projection, making the relationship between the parameter and the coefficients smoother compared to the linear counterpart (POD), as shown earlier in Figures \ref{fig:kpod_coeffs} and \ref{fig:pod_coeffs}.
This enhanced smoothness facilitates the use of simple yet effective regression or interpolation methods, such as polynomial regression of an appropriate degree or cubic spline interpolation (see Figure \ref{fig:polynomial_vs_spline}).

Finally, it is sufficient to compose the regressed parameter-to-coefficient map, with the decoding map \(\psi^D\) (see Section \ref{sec:ae-based-approach}), to approximate solutions to the problem for different values of the parameter. The composition \(\psi^D \circ \mathscr{F}_{\text{online}}: \mu \rightarrow u(\mu)\) maps parameters directly to the solution space, effectively serving as a complete solution predictor to Problem \eqref{eq:burgers-inviscid}.
We investigated several mappings corresponding to different polynomial degrees in terms of their computational efficiency and accuracy in representing the behavior of higher-order kPOD components. The optimal degree, equals to 25, was determined through an iterative process tailored to the specific problem, aimed at ensuring a balance between accuracy and model complexity to prevent overfitting, with performance primarily assessed by reconstruction fidelity. It's important to note, however, that a lower degree could have been used, which would have only resulted in a slight increase in reconstruction error.
The mean test error using a polynomial regression of degree 25 for $\mathscr{F}_{\text{online}}$ is \(\overline{\varepsilon}_{\text{online}} = 1.2\%\), which aligns perfectly with the compression-related reconstruction error \(\overline{\varepsilon}_{\text{kPOD}}\). This indicates that in this case the regression task on the coefficients does not introduce additional sources of error into the model.
The relative error fields, shown in Figure \ref{fig:burgers_inviscid}, visually illustrate deviations across a subset of tested parameters. While the autoencoder effectively reconstructs regions away from the shock, such as the rarefaction wave, minor misalignments near the discontinuity can result in significantly higher localized errors, particularly evident along the shock line.
\section{Conclusions and perspectives}%
\label{sec:conclusions_and_perspectives}

In this work, we introduced a novel reduced-order modeling framework that combines kernel POD, deep learning, and optimal transport theory to solve high-dimensional parametric PDEs. Our primary goal was to provide an efficient and accurate reduced-order representation strategy for solutions across a wide range of parameter values.

The foundation of our approach is the extension of traditional POD through the utilization of a kernel based on the Wasserstein distance. This kernel is efficiently approximated through entropic regularization and numerically optimized using the Sinkhorn algorithm. This way, the kPOD retains essential nonlinear features that might be lost by linear dimensionality reduction techniques.

To construct the inverse map of kPOD, we employed an autoencoder trained to learn the kPOD reduction. Following the autoencoder's training, it can be used to reconstruct the original solutions from the reduced coefficients. To ensure stability and convergence during the training process, we utilized the Sinkhorn divergence, which proved beneficial for training the proposed architectures, particularly showing improvements over the classical MSE in the case of convolutional layers.

We tested the proposed architecture on several benchmarks affected by a slow-decay of the Kolmogorov n-width, including the Poisson equation, the advection-diffusion problem, and the Burgers equation, revealing the framework's ability to obtain accurate and efficient reduced representations of solutions while capturing the geometric characteristics of the data.

The promising results obtained, pave the way for many interesting investigations and future improvements. Exploring alternative kernel functions and evaluating their performance for different PDE problems could provide additional insights into capturing specific features \cite{kernel-methods}.

Another interesting development could be in the direction of noisy data, adapting the framework to a variational autoencoder architecture capable of providing statistical interpretation, and characterized by a more robust training phase \cite{goodfellow}.

The assumption on the structured mesh is essential to guarantee the locality of the convolutional layers, which is lost for complex domain defined over unstructured grids. A further step towards a geometrically consistent and interpretable framework could originate by investigating convolutional autoencoders defined over graph neural networks \cite{PichiGraphConvolutionalAutoencoder2023}.

In addition, employing the kPOD-AE approach to create reduced-order models for parametric optimization and control problems has the potential to significantly accelerate optimization and control processes.

A critical avenue for future research lies in addressing the data-intensive nature of autoencoder training for high-dimensional PDE solutions. Exploring data-efficient learning techniques, such as transfer learning or active learning, holds promise for mitigating the demands of training data and enhancing the effectiveness of reduced-order modeling strategies \cite{transfer-learning,active-learning}.

One possibility lies in leveraging clustering techniques based on Wasserstein barycenters \cite{mula2022,mula2023}, to extract representative measures from the data. These extracted measures could then be employed as inputs for the kPOD approach.

Another strategy involves fitting the probability distributions using Gaussian models, in order to use known optimal transport maps \cite{taddei2022}. In particular, after training the model, it can be used to enrich the dataset used to construct the reduced-order basis for a projection-based approach \cite{iollo-augmentation}.

Finally, generalizing the framework to handle more complex physics and geometries would further improve its applicability to a broader range of problems in engineering and science.

\section*{Acknowledgment}
The authors would like to extend sincere appreciation to Prof. Jan S. Hestaven from \'Ecole Polytechnique F\'ed\'erale de Lausanne for generously hosting the visit of M. Khamlich at his institution. The collaborative opportunity during this visit significantly contributed to the successful realization of this project.

This work was partially funded by European Union Funding for Research and Innovation — Horizon 2020 Program — in the framework of European Research Council Executive Agency: H2020 ERC CoG 2015 AROMA-CFD project 681447 "Advanced Reduced Order Methods with Applications in Computational Fluid Dynamics" P.I. Professor Gianluigi Rozza. We also acknowledge the PRIN 2017 "Numerical Analysis for Full and Reduced Order Methods for the efficient and accurate solution of complex systems governed by Partial Differential Equations" (NA-FROM-PDEs).

 \bibliographystyle{habbrv}
\bibliography{references}
\newpage
\appendix
\section{NN architectures}%
\label{app:arch}
This section summarizes the neural network architectures used for the kPOD backward map approximation.
In Tables \ref{arch:conv} and \ref{arch:ff}, we outline the architectures for the convolutional and the feedforward autoencoder, respectively.
We recall that the input and output dimensions for a convolutional layer can be calculated using the following formulas:
\[
H^{l+1} = \frac{{H^{l} - F + 2P_0}}{S_0} + 1 \quad \text{and} \quad W^{l+1} = \frac{{W^{l} - F + 2P_1}}{S_1} + 1
\]
where $(H^{l},W^{l})$ and $(H^{l+1},W^{l+1})$ are the input and output dimensions, $F$ is the filter size, $S = (S_0, S_1)$ is the stride, and $P = (P_0, P_1)$ is the padding.

Similarly, the input and output dimensions for a transpose convolutional (deconvolutional) layer can be calculated using the following formulas:
\[
H^{l+1} = S_0 \times (H^{l} - 1) + F - 2P_0 \quad \text{and} \quad W^{l+1} = S_1 \times (W^{l} - 1) + F - 2P_1.
\]
In particular, the input to the first convolutional layer is such that $H^{1}\cdot W^{1}=N_{h}$ which is the number of dofs for the problem at hand.

In the feedforward architecture, dropout regularization is employed to avoid overfitting and enhance the model's generalization capability.
\begin{table}[ht]
\caption{Convolutional Autoencoder Architecture}\label{table:autoencoder-architecture}
\begin{tabular}{ccccccccc}
\hline
\rowcolor[gray]{0.8} Module & Layer & Input Size & Output Size & Activation & Stride & Padding \\
\hline
\multirow{5}{*}{Encoder} & \cellcolor[gray]{0.9}$\text{Conv1}$ &\cellcolor[gray]{0.9}\cellcolor[gray]{0.9}$n_b \times 1 \times H^1 \times W^1$ &\cellcolor[gray]{0.9}\cellcolor[gray]{0.9}$n_b \times 8 \times H^2 \times W^2$ &\cellcolor[gray]{0.9}\cellcolor[gray]{0.9}\cmark &\cellcolor[gray]{0.9}\cellcolor[gray]{0.9}$2 \times 2$ &\cellcolor[gray]{0.9}\cellcolor[gray]{0.9}$1 \times 1$ \\
& $\text{Conv2}$ & $n_b \times 8 \times H^2 \times W^2$ & $n_b \times 16 \times H^3 \times W^3$ & \cmark & $2 \times 2$ & $1 \times 1$ \\
& \cellcolor[gray]{0.9}$\text{Conv3}$ &\cellcolor[gray]{0.9}$n_b \times 16 \times H^3 \times W^3$ &\cellcolor[gray]{0.9}$n_b \times 32 \times H^4 \times W^4$ &\cellcolor[gray]{0.9}\cmark &\cellcolor[gray]{0.9}$2 \times 2$ &\cellcolor[gray]{0.9}$1 \times 1$ \\
& $\text{Conv4}$ & $n_b \times 32 \times H^4 \times W^4$ & $n_b \times 64 \times H^5 \times W^5$ & \cmark & $2 \times 2$ & $1 \times 1$ \\
& \cellcolor[gray]{0.9}$\text{Flatten}$ &\cellcolor[gray]{0.9}$n_b \times 64 \times H^5 \times W^5$ &\cellcolor[gray]{0.9}$n_b \times (64 \cdot H^5 \cdot W^5)$ &\cellcolor[gray]{0.9}\xmark &\cellcolor[gray]{0.9}- &\cellcolor[gray]{0.9}- \\
& $\text{FC1}$ & $n_b \times (64 \cdot H^5 \cdot W^5)$ & $n_b \times \text{k}$ & \cmark & - & - \\
\hline
\multirow{7}{*}{Decoder} & \cellcolor[gray]{0.9}$\text{FC\_out1}$ &\cellcolor[gray]{0.9}$n_b \times \text{k}$ &\cellcolor[gray]{0.9}$n_b \times (256)$ &\cellcolor[gray]{0.9}\cmark &\cellcolor[gray]{0.9}- &\cellcolor[gray]{0.9}- \\
& $\text{FC\_out2}$ & $n_b \times 256$ & $n_b \times (64 \cdot H^6 \cdot W^6 )$ & \cmark & - & - \\
& \cellcolor[gray]{0.9}$\text{Unflatten}$ &\cellcolor[gray]{0.9}$n_b \times (64 \cdot H^6 \cdot W^6)  $ &\cellcolor[gray]{0.9}$n_b \times 64 \times H^6 \times W^6$ &\cellcolor[gray]{0.9}\xmark &\cellcolor[gray]{0.9}- &\cellcolor[gray]{0.9}- \\
& $\text{ConvT1}$ & $n_b \times 64 \times H^6 \times W^6$ & $n_b \times 32 \times H^7 \times W^7$ & \cmark & $2 \times 2$ & $1 \times 1$ \\
& \cellcolor[gray]{0.9}$\text{ConvT2}$ &\cellcolor[gray]{0.9}$n_b \times 32 \times H^7 \times W^7$ &\cellcolor[gray]{0.9}$n_b \times 16 \times H^8 \times W^8$ &\cellcolor[gray]{0.9}\cmark &\cellcolor[gray]{0.9}$2 \times 2$ &\cellcolor[gray]{0.9}$1 \times 1$ \\
& $\text{ConvT3}$ & $n_b \times 16 \times H^8 \times W^8$ & $n_b \times 8 \times H^9 \times W^9$ & \cmark & $2 \times 2$ & $1 \times 1$ \\
& \cellcolor[gray]{0.9}$\text{ConvT4}$ &\cellcolor[gray]{0.9}$n_b \times 8 \times H^{10} \times W^{10}$ &\cellcolor[gray]{0.9}$n_b \times 1 \times H^1 \times W^1$ &\cellcolor[gray]{0.9}\xmark &\cellcolor[gray]{0.9}$1 \times 1$ &\cellcolor[gray]{0.9}$1 \times 1$ \\
\hline
\end{tabular}
\label{arch:conv}
\end{table}

\begin{table}[ht]
\caption{Fully Connected Feedforward Autoencoder Architecture}\label{table:feedforward-autoencoder-architecture}
\begin{tabular}{ccccc}
\hline
\rowcolor[gray]{0.8}Module & Layer & Input Size & Output Size & Activation \\
\hline
\multirow{7}{*}{Encoder}&\cellcolor[gray]{0.9}$\text{FC1}$ &\cellcolor[gray]{0.9} $n_b \times N_{h}$ &\cellcolor[gray]{0.9} $n_b \times 64$ &\cellcolor[gray]{0.9}\cmark \\
               &$\text{FC2}$ & $n_b \times 64$ & $n_b \times 128$ & \cmark \\
               &\cellcolor[gray]{0.9}$\text{FC3}$ &\cellcolor[gray]{0.9}$n_b \times 128$ &\cellcolor[gray]{0.9}$n_b \times 256$ &\cellcolor[gray]{0.9}\cmark \\
               &$\text{FC4}$ & $n_b \times 256$ & $n_b \times 256$ & \cmark \\
               &\cellcolor[gray]{0.9}$\text{FC5}$ &\cellcolor[gray]{0.9}$n_b \times 256$ &\cellcolor[gray]{0.9}$n_b \times 128$ &\cellcolor[gray]{0.9}\cmark \\
               &$\text{FC6}$ & $n_b \times 128$ & $n_b \times 64$ & \xmark \\
               &\cellcolor[gray]{0.9}$\text{FC7}$ &\cellcolor[gray]{0.9}$n_b \times 64$ &\cellcolor[gray]{0.9}$n_b \times k$ &\cellcolor[gray]{0.9}\xmark \\
\hline
 \multirow{7}{*}{Decoder}&$\text{FC\_out1}$ & $n_b \times k$ & $n_b \times 64$ & \cmark \\
             &\cellcolor[gray]{0.9}$\text{FC\_out2}$ &\cellcolor[gray]{0.9}$n_b \times 64$ &\cellcolor[gray]{0.9}$n_b \times 128$ &\cellcolor[gray]{0.9}\cmark\\
             &$\text{FC\_out3}$ & $n_b \times 128$ & $n_b \times 256$ & \cmark \\
             &\cellcolor[gray]{0.9}$\text{FC4\_out4}$ &\cellcolor[gray]{0.9}$n_b \times 256$ &\cellcolor[gray]{0.9}$n_b \times 256$ &\cellcolor[gray]{0.9}\cmark \\
             &$\text{FC5\_out5}$ & $n_b \times 256$ & $n_b \times 128$ & \cmark \\
             &\cellcolor[gray]{0.9}$\text{FC\_out6}$ &\cellcolor[gray]{0.9}$n_b \times 128$ &\cellcolor[gray]{0.9}$n_b \times 64$ &\cellcolor[gray]{0.9}\xmark \\
             &$\text{FC\_out7}$ & $n_b \times 64$ & $n_b \times N_{h}$ & \xmark \\
\hline
\end{tabular}
\label{arch:ff}
\end{table}

\end{document}